\documentclass[11pt]{article}

\usepackage[top=1in,bottom=1.2in,left=1.25in,right=1.25in]{geometry}

\usepackage{graphicx}
\usepackage{ragged2e}
\usepackage{tabto}
\usepackage{amssymb}
\usepackage{float}
\usepackage{multirow}
\usepackage{amsmath}
\usepackage{tikz}
\usepackage{amsthm}
\usepackage{color}
\usepackage{url}
\usepackage[noabbrev,capitalise]{cleveref}
\usepackage{afterpage}
\usepackage{multicol}
\usepackage[font=footnotesize]{caption}
\usepackage[font=footnotesize]{subcaption}
\usepackage{marginnote}
\usepackage{tcolorbox}

\usepackage{newtxtext}
\usepackage{newtxmath}

\usepackage[absolute,overlay]{textpos}

\usepackage[sf,bf,medium]{titlesec}

\usepackage{bm}
\usepackage{dsfont}
\usepackage{booktabs}
\usepackage{siunitx}

\newcommand{\lp}{\left(}
\newcommand{\rp}{\right)}

\usepackage[sort,compress]{cite}

\newcommand\bx{\boldsymbol{x}}

\newcommand\xu{\partial_{u} \boldsymbol{x}}
\newcommand\xv{\partial_{v} \boldsymbol{x}}
\newcommand\bs{\boldsymbol{s}}
\newcommand\by{\boldsymbol{y}}

\newcommand\bH{\boldsymbol H}
\newcommand\bF{\boldsymbol F}

\newcommand\bn{\boldsymbol n}

\newcommand\bzeta{\boldsymbol \zeta}
\newcommand\bmu{\boldsymbol \mu}

\newcommand\bX{\boldsymbol X}

\newcommand\bV{\boldsymbol V}

\newcommand\bv{\boldsymbol v}
\newcommand\bt{\boldsymbol t}

\newcommand\surflap{\Delta_\Gamma}

\newcommand\bR{\mathbb{R}}
\newcommand\niter{n_{\textrm{iter}}}

\newcommand\Npat{N_{\textrm{patches}}}

\newcommand\epsgmres{\varepsilon_{\textrm{GMRES}}}

\newcommand\Nover{N_{\textnormal{over}}}
\newcommand\wover{w^{\textnormal{over}}}

\newtheorem{remark}{\sffamily Remark}

\newtheorem{lemma}{\sffamily Lemma}

\newcommand{\cW}{\mathcal W}
\newcommand{\cD}{\mathcal D}
\newcommand{\cS}{\mathcal S}
\newcommand{\cA}{\mathcal A}
\newcommand{\cK}{\mathcal K}

\newcommand{\cI}{\mathcal I}

\newcommand{\ucomp}{u_{\textrm{comp}}}

\numberwithin{equation}{section}

\begin{document}

%\maketitle

\begin{titlepage}

  \raggedleft
  %{\sffamily \bfseries STATUS: {\color{red} Unsubmitted manuscript}}
  {\sffamily \bfseries STATUS: arXiv pre-print}
  
  \hrulefill

  \raggedright
  \begin{textblock*}{\linewidth}(1.25in,2in) % {block width} (coords) 
    {\LARGE \sffamily \bfseries FMM-accelerated solvers for
      the Laplace-Beltrami\\
      \vspace{.25\baselineskip}
      problem on complex surfaces in three dimensions}
  \end{textblock*}

  \normalsize

  \vspace{2in}
  Dhwanit Agarwal\\
  \emph{\small Oden Institute, University of Texas at Austin\\
    Austin, TX, 78712}\\
  \texttt{\small dhwanit@oden.utexas.edu}

   \vspace{\baselineskip}
  Michael O'Neil\\
  \emph{\small Courant Institute, NYU\\
    New York, NY 10012}\\
  \texttt{\small oneil@cims.nyu.edu}
  
  \vspace{\baselineskip}
  Manas Rachh\\
  %\footnote{Corresponding author.}\\
  \emph{\small Center for Computational Mathematics, Flatiron Institute\\
    New York, NY 10010}\\
  \texttt{\small mrachh@flatironinstitute.org}

  \begin{textblock*}{\linewidth}(1.25in,7in) % {block width} (coords) 
    \today
  \end{textblock*}

\end{titlepage}

\begin{abstract}
  The Laplace-Beltrami problem on closed surfaces embedded in three dimensions
  arises in many areas of physics, including molecular dynamics (surface
  diffusion), electromagnetics (harmonic vector fields), and fluid dynamics
  (vesicle deformation). Using classical potential theory, the Laplace-Beltrami
  operator can be pre-/post-conditioned with an integral operator whose kernel
  is translation invariant, resulting in well-conditioned Fredholm integral
  equations of the second-kind. These equations have the standard~$1/ r$ kernel
  from potential theory, and therefore the equations can be solved rapidly and
  accurately using a combination of fast multipole methods (FMMs) and high-order
  quadrature corrections.  In this work we detail such a scheme, presenting two
  alternative integral formulations of the Laplace-Beltrami problem, each of
  whose solution can be obtained via FMM acceleration. We then present several
  applications of the solvers, focusing on the computation of what are known as
  \emph{harmonic vector fields}, relevant for many applications in
  electromagnetics. A battery of numerical results are presented for each
  application, detailing the performance of the solver in various geometries. \\
  

\end{abstract}

\tableofcontents

\newpage

\section{Introduction}

The Laplace-Beltrami problem along surfaces in three dimensions appears in many
fields, including plasma physics~\cite{beltrami3}, fluid
mechanics~\cite{rahimian2010petascale,veerapaneni_2011},
electromagnetics~\cite{epstein-2010,gendeb2,superconductor1,nedelec}, surface
diffusion, pattern formation~\cite{kim2020pattern,plaza2004effect}, and
computational geometry~\cite{kromer2018highly}. In general, the Laplace-Beltrami
operator is a second-order variable coefficient elliptic operator whose behavior
is highly dependent on the nature of the surface on which it is defined. Its
eigenfunctions are useful in describing band-limited functions defined along the
associated manifold or surface, and more generally, the related harmonic vector
fields provide a basis for important classes of functions appearing in
electromagnetics (e.g. non-radiating currents, etc.).

To describe the problem, suppose that $\Gamma$ is at least a
twice-differentiable surface without boundary embedded in~$\mathbb{R}^{3}$. Let
$\nabla_{\Gamma}$ denote the intrinsic surface gradient, and let its adjoint
acting on the space of continuously differentiable $\mathbb{L}^2$ tangential
vector fields on the boundary, denoted by $\nabla_{\Gamma} \cdot$, be the
intrinsic surface divergence. Given a mean-zero function $f$ defined
along~$\Gamma$, the Laplace-Beltrami problem is to find a mean-zero function $u$
satisfying $\Delta_{\Gamma} u = f$, where $\Delta_{\Gamma}$ is the
Laplace-Beltrami operator defined as $\Delta_{\Gamma}= \nabla_{\Gamma} \cdot
\nabla_{\Gamma}$. The mean-zero condition on~$f$ and $u$ is required for unique
solvability due to a rank-deficiency of the Laplace-Beltrami
operator~\cite{oneil2018,epstein-2010,nedelec}.

In many of the applications discussed above, the Laplace-Beltrami problem
naturally arises for the computation of the Hodge decomposition of tangential
vector fields. Any smooth tangential vector field $\bF$ on the surface $\Gamma$
can be written as
\begin{equation}
\bF = \nabla_{\Gamma} \alpha + \bn \times \nabla_{\Gamma} \beta + \bH \, ,
\end{equation}
where $\alpha$, $\beta$ are smooth functions on $\Gamma$, $\bH$ is a harmonic
field satisfying $\nabla_{\Gamma} \cdot \bH = \nabla_{\Gamma} \cdot \bn \times
\bH = 0$, and $\bn$ is the unit normal to the surface $\Gamma$. The curl free
component $\nabla_{\Gamma} \alpha$ and the divergence free component $\bn \times
\nabla_{\Gamma} \beta$ can be computed via the solution of the following two
Laplace-Beltrami problems:
\begin{equation}
  \begin{aligned}
\Delta_{\Gamma} \alpha &= \nabla_{\Gamma} \cdot \bF, \\
\Delta_{\Gamma} \beta &= - \nabla_{\Gamma} \cdot (\bn \times \bF) .
\end{aligned}
\end{equation}
A related problem is that of determining a basis for harmonic vector fields on
non-contractible surfaces, which plays a crucial role for several integral
representations in computational electromagnetics. On a genus $g$ object, the
space of harmonic vector fields forms a $2g$-dimensional vector space. One can
compute elements of this basis via the solution of the two Laplace-Beltrami
problems above for a randomly generated smooth vector field $\bF$, which with
high probability tends to have a non-trivial projection onto the subspace of
harmonic vector fields.

Numerical methods for the Laplace-Beltrami problem can be typically classified
into two categories.  The first category of methods rely on the direct
discretization of the differential operator on surface. Examples include finite
element
methods~\cite{bansch2005finite,bonito2020finite,burman2017cut,demlow2007adaptive},
virtual element methods~\cite{beirao2014hitchhiker,frittelli2018virtual},
differencing methods~\cite{wang2018modified}, and level set
methods~\cite{bertalmio2001variational,chen2015closest,greer2006fourth,macdonald2008level,macdonald2010implicit}.
The level set methods differ from the rest of the direct discretizations of the
PDE as they rely on using an embedded finite difference grid in the
volume~$\mathbb{R}^3$ to discretize a thickening of the surface.

On the other hand, the other general category of methods is to reformulate the
Laplace-Beltrami problem as an integral equation along the surface. For example,
on a flat surface the Laplace-Beltrami operator reduces to the standard Laplace
operator and therefore the Green's function is known to be~$G(\bx) =
-\log{|\bx|}/(2\pi)$. This Green's function~$G$ can be used as a parametrix for
the fundamental solution to the Laplace-Beltrami operator on surfaces as
well~\cite{kropinski_2014_fast}. This idea has also been used for the solution
of the related Yukawa-Beltrami problem on general surfaces, including surfaces
with boundary in~\cite{kropinski_2016_integral}. A separation of variables
approach, combined with an integral reformulation of an ODE, is presented
in~\cite{goodwill2021numerical}; this approach was also shown to be applicable
in geometries with edges. An alternate integral equation-based approach is to
use the Green's function of the free-space Laplace equation in three dimensions
as either left or right preconditioners for the Laplace-Beltrami
problem~\cite{oneil2018}.

In this work, we present a fast multipole method (FMM) accelerated iterative
solver for the solution of the Laplace-Beltrami problems on complex
three-dimensional surfaces based on the integral formulation presented
in~\cite{oneil2018}, and using the locally corrected quadrature approach
of~\cite{greengard2020fmm}. In some respects, this work can be considered
\emph{Part II} of the earlier work~\cite{oneil2018}. The primary motivation of
choosing this approach is two-fold: (1) the resulting integral equations tend to
be as well-conditioned as the underlying physical problem (which is quite
well-conditioned in this case, in an absolute sense), and (2) since the
representation uses layer potentials related to the free space Green's function
in three dimensions, the numerical solver can be coupled to existing optimized
FMM libraries~\cite{greengard-1997, greengard2020fmm}. We demonstrate the
well-conditioned nature of the integral equation through a battery of numerical
experiments, each of which incorporates a high-order representation of the
surface and high-order quadratures. We also show the importance of what we refer
to as \emph{numerical second-kindness} when discretizing the integral equation
using locally corrected quadrature methods. These methods rely on the accurate
computation of the principal-value part of the layer potentials in the
representation. In particular, while there are several equivalent integral
formulations of the same PDE on the surface, we show that the formulation which
explicitly handles the identity term of the integral equation tends to be the
most stable (from a numerical point of view). Lastly, with the acceleration
provided by FMMs, our integral formulation can be used to compute a basis for
harmonic vector fields in relatively high-genus geometries.

The rest of the paper is organized in follows. In~\cref{sec:lapbel}, we describe
the Laplace-Beltrami problem.  In~\cref{sec:inteq}, we review the integral
equation formulation presented in~\cite{oneil2018}, and in~\cref{sec:solver}, we
review a high order FMM-accelerated iterative solver. In~\cref{sec:harm-vec}, we
discuss an application of the solver to the computation of a basis for harmonic
vector fields on non-contractible surfaces. In~\cref{sec:num}, we present
several numerical examples to demonstrate the efficiency of our approach,
including a discussion on the importance of numerical second-kindness. Finally,
in~\cref{sec:conclusions}, we conclude by discussing avenues for future
research.

\section{The Laplace-Beltrami problem \label{sec:lapbel}}
The Laplace-Beltrami operator, also known as the surface Laplacian, is
the generalization of the Laplace operator to general (smooth)
manifolds.  In this work, we will focus on smooth surfaces embedded in
three dimensions.  Suppose $\Gamma \subset \mathbb{R}^{3}$ is an
oriented surface with at least two continuous derivatives, and let
$\bx(u,v): D \subset \mathbb{R}^{2} \to \mathbb{R}^{3}$ denote a local
parametrization of~$\Gamma$. The first fundamental form, or
equivalently the metric tensor at $\bx(u,v)$, is given by
\begin{equation}
    g = \begin{bmatrix}
    g_{11}       & g_{12}  \\
    g_{21}       & g_{22} 
\end{bmatrix} 
= \begin{bmatrix}
    \xu \cdot \xu       &  \xu \cdot \xv \\
    \xv \cdot \xu     & \xv \cdot \xv
\end{bmatrix} ,
\end{equation}
where $\xu$ and $\xv$ are abbreviations for $\frac{\partial \bx}{\partial u}$
and $\frac{\partial \bx}{\partial v}$, respectively, and are tangent vectors to the
surface at $\bx \in \Gamma$. The unit normal at $\bx$ can written as
\begin{equation}
  \bn(\bx) = \frac{\xu \times \xv}{ | \xu \times \xv |}.
\end{equation}
The unit vectors in the directions of $\xu$, $\xv$, and $\bn(\bx)$ form
a local coordinate system (not necessarily orthogonal) at each
point~$\bx\in \Gamma$.
 
Next, suppose that $\boldsymbol{F}$ is a smooth tangential vector
field defined along~$\Gamma$ and is expressed in terms of the local
tangent vectors $\xu$, $\xv$ as
 \begin{equation}
 \boldsymbol{F}(u,v) = F_{1} \xu + F_{2} \xv \, .
 \end{equation} 
The surface divergence of $\boldsymbol{F}$, denoted by
$\nabla_{\Gamma} \cdot \boldsymbol{F}$, is then~\cite{nedelec,frankel} given by
 \begin{equation}
 \nabla_{\Gamma} \cdot \boldsymbol{F} = \frac{1}{\sqrt{|g|}} \left( 
 \partial_{u} \left( \sqrt{|g|} F_{1} \right) + \partial_{v} \left( \sqrt{|g|} F_{2} \right)
 \right)  \, ,
 \end{equation}
 where $|g|$ is the determinant of the metric tensor $g$.  The surface
 gradient operator $\nabla_{\Gamma}$, formally defined as the adjoint of the surface divergence operator on the space
 of~$\mathbb{L}^{2}$ functions with respect to the induced metric on $\Gamma$,
 is given by
 \begin{equation}
 \nabla_{\Gamma} \psi = \left( g^{11} \partial_{u} \psi + g^{12} \partial_{v} \psi \right) \xu + 
 \left( g^{21} \partial_{u} \psi + g^{22} \partial_{v} \psi \right) \xv \, , 
 \end{equation}
 where $\psi = \psi(u,v)$ is a scalar function defined on~$\Gamma$
 and $g^{ij}$ denote the components of the inverse of the
 metric tensor, denoted by $g^{-1}$.
 
The Laplace-Beltrami operator is then defined as~$\surflap =
\nabla_{\Gamma} \cdot \nabla_{\Gamma}$, and is explicitly given by
the expression
\begin{equation}
     \Delta_{\Gamma} \psi = \frac{1}{\sqrt{|g|}} \left( \partial_{u}
     \left[ \sqrt{|g|} \left(g^{11} \partial_{u} \psi + g^{12}
       \partial_{v} \psi \right) \right] + \partial_{v} \left[
       \sqrt{|g|} \left(g^{21} \partial_{u} \psi + g^{22} \partial_{v}
       \psi \right) \right] \right).
\end{equation}
The standard Laplace-Beltrami problem is to solve the following
PDE for~$\psi$ along~$\Gamma$,
\begin{equation}
  \label{eq:lapbel}
  \surflap \psi = f,
\end{equation}
where~$f$ is a given function. Conditions on the regularity required
in~$f$, and the subsequent regularity obtained in~$\psi$ have been
studied in~\cite{rosenberg, jost,epstein-2010, goodwill2021numerical}.
We will assume that all functions and surfaces are
sufficiently smooth so as to obtain high order accuracy in the
subsequently numerical examples.

As is clear, the Laplace-Beltrami operator maps constant functions
along~$\Gamma$ to zero. Thus, the equation~\eqref{eq:lapbel} is
rank-one deficient~\cite{imbertgerard_2017, sifuentes_2015} (and
formally self-adjoint). The
following lemma~\cite{rosenberg, oneil2018} summarizes a classical
result regarding the well-posedness of the above Laplace-Beltrami
problem.

\begin{lemma}
 \label{lem:meanzero}
Let $\Gamma$ be a smooth, closed, and orientable
boundary. Along~$\Gamma$, the
Laplace-Beltrami operator is uniquely invertible as a map from
$\mathbb{M}_{0} \to \mathbb{M}_{0}$, the space of mean-zero functions
defined on $\Gamma$, \,
\begin{equation}
\mathbb{M}_{0} = \left\{ \psi \, : \, \int_{\Gamma} \psi \, da = 0 \right\} \, .
\end{equation}
Indeed, if the function~$f$ defined along~$\Gamma$ has mean zero, then
there exists a unique mean-zero solution,~$\psi$, to the
Laplace-Beltrami problem~$\surflap \psi = f$.
\end{lemma}

\section{Integral equation formulations}
\label{sec:inteq}

In this section, we review an equivalent well-conditioned second-kind integral
equation for the solution of the Laplace-Beltrami problem. This formulation
relies on preconditioning the Laplace-Beltrami operator on both the left and
right using the Laplace single layer potential. The resulting Fredholm operator
can be expressed in terms of standard layer potentials for the free space
Laplace equation in three dimensions. Before describing this formulation, we
first define the relevant layer potentials used in the formulation. The Green's
function $G$ for Laplace's equation in three dimensions is given by
\begin{equation}
G(\bx,\by) = \frac{1}{4\pi |\bx-\by| }  \, .
\end{equation}
The  single and double layer potentials are then given by:
\begin{equation}
  \begin{aligned}
  \mathcal{S}[\sigma](\bx) &= \int_{\Gamma} G(\bx,\by) \, 
      \sigma(\by) \, da, \\
  \mathcal{D}[\mu](\bx) &= \int_{\Gamma} \lp \nabla_{\by} G(\bx,\by)
\cdot \bn(\by) \rp
\mu(\by) \, da.
  \end{aligned}
\end{equation}
When the \emph{target} point $\bx$ lies on~$\Gamma$, the integral
defining the double layer should be
interpreted in a principal value sense. For smooth surfaces $\Gamma$,
both $\mathcal{S}$ and $\mathcal{D}$ are compact operators as maps
from~$\mathbb{L}^{2}(\Gamma) \to \mathbb{L}^{2}(\Gamma)$.
 Related to the single and double layer
potentials are the restrictions of their normal derivatives to
$\Gamma$ given by
\begin{equation}
  \begin{aligned}
    \mathcal{S}'[\sigma](\bx) &= \bn(\bx) \cdot \nabla \int_{\Gamma}
    G(\bx,\by) \sigma(\by) \, da, \\
    \mathcal{D}'[\mu](\bx) &= \int_{\Gamma} \bn(\bx) \cdot
    \nabla_{\bx} \nabla_{\by}  G(\bx,\by) \cdot \bn(\by) \mu(\by)
    \, da.
  \end{aligned}
\end{equation}
The integral defining $\mathcal{S}'$ is to be interpreted in a
principal value sense, and that for $\mathcal{D}'$ is to be
interpreted in a finite-part sense~\cite{kress_2014, nedelec}. The
operator $\mathcal{S}'$ is also a compact operator on
$\mathbb{L}^{2}(\Gamma)$, whereas~$\mathcal D'$ is a hypersingular
operator.

Finally, let $\mathcal{S}''$ denote the restriction of the second
normal derivative of $\mathcal{S}$ to $\Gamma$ given by:
\begin{equation}
  \mathcal{S}''[\sigma](\bx) = \int_{\Gamma} \lp
  \bn(\bx) \cdot \nabla_{\bx} \nabla_{\bx}  G(\bx,\by) \cdot \bn(\bx)
  \rp \sigma(\by) \, da,
\end{equation}
where, as for the integral defining $\mathcal{D}'$, the integral
should be interpreted in a finite-part sense \cite{kress_2014,
  nedelec, oneil2018}. It can be shown that the \emph{difference}
operator $\mathcal{S}'' + \mathcal{D}'$ is a compact operator
on~$\mathbb{L}^{2}(\Gamma)$~\cite{nedelec}, as the dominant
singularities in the integrand cancel each other out.

We now review two equivalent second-kind boundary integral formulations for the
solution of the Laplace-Beltrami equation \eqref{eq:lapbel} as presented
in~\cite{oneil2018}. One formulation is to represent the solution $\psi =
\cS[\sigma]$, where $\sigma$ is an unknown density and precondition to
equation~\eqref{eq:lapbel} using $\cS$, i.e. solve
\begin{equation}
\cS \Delta_{\Gamma} \cS [\sigma] = \cS[f] \,.
\end{equation}
Using Calder\'on identities, the equation above can be rewritten as
\begin{equation}
  \label{eq:bie0}
-\frac{\sigma}{4} + \mathcal{D}^2[\sigma] -
\mathcal{S}\left(\mathcal{S}'' + \mathcal{D}' + 2H\mathcal{S}'
\right)[\sigma] = \cS[f] \, ,
\end{equation}
where $H$ is the mean curvature along~$\Gamma$. The above equation still has a
one dimensional null-space which is related to constant functions comprising the
null-space of the Laplace Beltrami problem. To impose the integral equation on
the space of mean-zero functions, we instead solve the integral equation
\begin{equation}
\label{eq:bie1}
-\frac{\sigma}{4} + \mathcal{D}^2[\sigma] -
\mathcal{S}\left(\mathcal{S}'' + \mathcal{D}' +2H\mathcal{S}'
\right)[\sigma] + \mathcal{SWS}[\sigma] = \mathcal{S}[f] \, ,
\end{equation}
where $\mathcal{W}$ is the operator given by
\begin{equation}
\mathcal{W}[\sigma] = \frac{1}{|\Gamma|} \int_{\Gamma} \sigma \, da \,.
\end{equation}
If $f$ is a mean-zero function in $\mathbb{L}^{2} (\Gamma)$, then there exists a unique
solution $\sigma$ in $\mathbb{L}^{2}(\Gamma)$ which satisfies~\cref{eq:bie1},  along with
$\int_{\Gamma} \cS[\sigma] \, da = 0$. Furthermore, $\psi = \cS[\sigma]$ is the unique mean-zero solution to
the Laplace-Beltrami problem $\Delta_{\Gamma} \psi = f$.

An alternate integral formulation is to instead represent $\psi = \cS^2 [\sigma]$, where
$\sigma$ as before is an unknown density and solve
\begin{equation}
\Delta_{\Gamma} \cS^2 [\sigma] = f \, ,
\end{equation}
which using Calder\'on identities, can be rewritten as
\begin{equation}
\label{eq:bie20}
-\frac{\sigma}{4} + (\cS')^2[\sigma] - (\cS''+\cD' + 2H\cS')\cS [\sigma] = f \, .
\end{equation}
The integral equation above has a similar one-dimensional null space issue as~\cref{eq:bie0}, which
can be addressed using the averaging operator $\cW$ by instead solving 
\begin{equation}
\label{eq:bie2}
-\frac{\sigma}{4} + (\cS')^2[\sigma] - (\cS''+\cD' + 2H\cS')\cS [\sigma] + \cW \cS^2 [\sigma]= f \, .
\end{equation}
Again, if $f$ is a mean-zero function in $\mathbb{L}^{2}(\Gamma)$, then there exists a unique solution
$\sigma$ in $\mathbb{L}^{2}(\Gamma)$ which satisfies~\cref{eq:bie2}, along with 
$\int_{\Gamma} \cS^2[\sigma] \, da = 0$. Furthermore, $\psi = \cS^2[\sigma]$ is the unique mean-zero solution
the Laplace-Beltrami problem $\Delta_{\Gamma} \psi = f$. 
A proof of these results can be found in~\cite{oneil2018}.

\section{FMM-accelerated solvers \label{sec:solver}}
\label{sec:fmmsolver}

In this section, we describe the details of a numerical solver for integral
equation formulations of the Laplace-Beltrami problem. For this, we first
briefly recap the FMM-accelerated locally corrected quadrature approach
of~\cite{greengard2020fmm}. For simplicity, consider a second-kind integral
equation of the form 
\begin{equation}
\label{eq:skie-demo}
\sigma(\bx) + \int_{\Gamma} K(\bx,\bx') \, \sigma(\bx') \, da(\bx') = f(\bx) \,
\quad \bx \in \Gamma \,,
\end{equation}
where $K(\bx,\bx')$ is either $G(\bx,\bx')$ or one its directional derivatives.
Suppose that the surface $\Gamma$ is represented via a disjoint union of patches
$\Gamma_{j}$ so that $\Gamma = \cup_{j=1}^{\Npat} \Gamma_{j}$, and that each
patch $\Gamma_{j}$ is parametrized by a non-degenerate chart $\bX^{j} : T_{0}
\to \Gamma_{j}$, where $T_{0}$ is the standard simplex given by 
\begin{equation*}
T_{0} = \{ (u,v): u \geq 0 \,, v\geq 0\, , u+v \leq 1 \} \subset \mathbb{R}^{2} \, .
\end{equation*}
An order $p$ discretization of $\Gamma$ is one where the components of $\bX^{j},
\partial_{u} \bX^{j}$, and $\partial_{v} \bX^{j}$ are represented  by a
polynomial expansion of total degree less than $p$ on $T_{0}$. One way to obtain
such a representation is to sample the chart and the derivative information at
order $p$ Vioreanu-Rokhlin nodes~\cite{vioreanu_2014} which are stable
interpolation nodes for orthogonal polynomials~\cite{bremer_2012,
koornwinder_1975} on $T_{0}$.

\begin{remark}
Triangulations of surfaces obtained from CAD or standard meshing packages often
tend to be low order and introduce several artificial edges and corners on the
discretized surface. In our examples, with the aim of showing high-order
convergence, we use the surface smoothing algorithm discussed
in~\cite{vico2020surface} to convert these low-order triangulations of complex
geometries to arbitrarily high order triangulations of a \emph{nearby}
related~$C^{\infty}$ surface. In some of our examples, analytic
parameterizations of the domains are known and are used instead of the algorithm
of~\cite{vico2020surface}. The algorithm was only applied to complex geometries,
such as the multi-holed object in Figure~\ref{fig:g10}.
\end{remark}

Let $ \{ \bx_{i}, w_{i} \}_{i=1}^{N}$ denote the union of the samples of the
charts $\bX^{j}$ at the Vioreanu-Rokhlin nodes of order $p$, and the
corresponding quadrature weights for integrating smooth functions on the
surface. A Nystr\"{o}m discretization of~\cref{eq:skie-demo} is given by 
\begin{equation}
\label{eq:nys}
\sigma_{i} + w_{ii} \sigma_{i} + \sum_{j\neq i} w_{ij} K(\bx_{i},\bx_{j}) \sigma_{j} = f(\bx_{i}) \, ,
\end{equation}
where $\sigma_{i}$ is an approximation for the solution $\sigma(\bx_{i})$, and
$w_{ij}$ are target dependent quadrature weights corresponding to a high-order
accurate discretization of the integral appearing in
equation~\eqref{eq:skie-demo}. In~\cite{greengard2020fmm}, for targets~$\bx_{i}
\in \Gamma_{\ell}$, this integral  is numerically approximated using a
combination of generalized Gaussian quadrature for the contribution
of~$\Gamma_{\ell}$ to the overall layer potential, adaptive integration for
patches $\Gamma_{j}$ which are \emph{close} to $\bx_{i}$, and target independent
oversampled quadratures for all the remaining patches (which all have smooth
integrands). For each discretization point $\bx_{i}$, we can split the
computational domain as
\begin{equation}
   \bR^{3} = \textrm{Far}(\bx_{i}) \cup \textrm{Near}(\bx_{i}),
\end{equation}
 with $\textrm{Near}(\bx_{i})$ containing $O(1)$ patches, and hence $O(1)$
discretization points. 

Let  $\{ \bs_{i}, \wover_{i} \}_{i=1}^{\Nover}$ denote the oversampled quadrature nodes and weights for integrating smooth functions on $\Gamma$, then the Nystr\"{o}m discretization in~\cref{eq:nys} takes the form
\begin{equation}
\label{eq:nys-lcq}
\sigma_i + \sum_{j=1}^{n_{i}} w_{ij} \, K(\bx_{i},
\by_{ij}) \, \tilde{\sigma}_{ij} + \sum_{\substack{j=1 \\ \bs_{j} \in
    \textrm{Far}(\bx_{i})}}^{\Nover} \wover_{j} \, K(\bx_{i}, \bs_{j}) \, 
\tilde{\sigma}_{j} = f(\bx_{i}) \, ,
\end{equation}
where~$\tilde{\sigma_{j}}$ and~$\tilde{\sigma}_{ij}$ correspond to (approximate)
interpolated values of~$\sigma$ evaluated at~$\bs_{j}$ and~$\by_{ij}$,
respectively, obtained via implicit local polynomial interpolation of the
samples~$\sigma_i$. The pairs~$\{\by_{ij}, w_{ij}\}$ are the quadrature nodes
and weights used in the adaptive integration procedure to accurately compute the
contribution of the layer potential due to the self panel~$\Gamma_{\ell}$ and
the associated near panels at the target location~$\bx_{i}$. After adding and
subtracting the contribution of the oversampled sources
in~$\textrm{Near}(\bx_{i})$, and then subsequently composing the interpolation
matrices to compute~$\tilde{\sigma}_{ij}$ and~$\tilde{\sigma}_{j}$
in~$\textrm{Near}(\bx_{i})$ and the corresponding kernel evaluations, the
discretized integral equation in~\eqref{eq:nys-lcq} can be rewritten as
\begin{equation}
  \label{eq:nys-lcq2}
  \sigma_i + \sum_{\substack{j=1\\\bx_{j} \in \textrm{Near}(\bx_{i})}}^{N} \tilde{w}_{ij} \, \sigma_j 
  + \sum_{j=1}^{\Nover} \wover_{j} \, K(\bx_{i}, \bs_{j}) \, 
\tilde{\sigma}_{j} = f(\bx_{i}),
\end{equation}
where~$\tilde{w}_{ij}$ represents an adjusted weight-kernel product obtained
from the addition/subtraction procedure. (See~\cite{greengard2020fmm} for a very
detailed discussion of this procedure, which we have very concisely summarized
above.)

The rate limiting step for numerically evaluating the layer potential via this
procedure  is the computation of the local quadrature corrections. However, when
using an an iterative algorithm such as GMRES~\cite{saad-1986} to solve the
integral equation, the effective quadrature-kernel weights~$\tilde{w}_{ij}$
above are precomputed and stored. Thus the computational cost of generating the
quadrature corrections is amortized over the number of GMRES iterations. The
memory requirement of this storage scales like $O(N)$ as there are $O(1)$ points
in $\textrm{Near}(\bx_{i})$. Since the far field part of the layer potential is
computed using a target independent quadrature rule, i.e. the sum over the
oversampled sources $\bs$ in~\eqref{eq:nys-lcq2}, the corresponding contribution
can be computed using standard fast multipole methods~\cite{greengard-1997,
greengard2020fmm}. Using this fast layer potential evaluator,  we can obtain the
solution~$\sigma$ in~$O(N)$ time if~\eqref{eq:nys-lcq} is well-conditioned. 

The recent release of the~\texttt{FMM3D} package allows for the computation of
the vector version of FMM sums, i.e., FMM sums with the same kernel, same source
and target locations, but with different strength vectors. Let $\bs_{j}$ denote
the source locations, $c_{\ell, j}$ denote the charge strengths, and
$\bv_{\ell, j}$ denote the dipole vectors, at $N$ sources, $j=1,2,\ldots N$, and~$\ell=1,2\ldots n_{d}$ densities, then the vector version of the Laplace FMM
computes the potentials $u_{\ell}(\bx)$ given by
\begin{equation}
  \label{eq:fmm-formula}
  u_{\ell}(\bx) 
  = \sum_{j=1}^{N} \frac{c_{\ell, j}}{|\bx - \bs_{j}|} - \bv_{\ell, j} \cdot \nabla_{\bx} \frac{1}{|\bx-\bs_{j}|} \, , \quad  \ell =1,2\ldots n_{d} \, ,
\end{equation}
along with the their gradients and hessians (if requested) at a given set of
target locations $x = \bt_{i}$, $i=1,2\ldots M$. The vector version of the FMM
is (empirically) computationally more efficient than separately calling $n_{d}$
separate FMMs since various density independent quantities can be reused in the
FMM algorithm. This feature of the FMM can be used for the more efficient
computation of multiple operators in the integral equations for the
Laplace-Beltrami problem.

We now turn our attention to discussing the specifics of GMRES accelerated
iterative solvers for the Laplace Beltrami problem. While piecewise smooth
discretizations of surfaces offer convenience to describe complicated three
dimensional surfaces, they pose an additional challenge when solving surface
PDEs such as the Laplace-Beltrami problem. In particular, the direct
discretization of the Laplace-Beltrami integral equation in the form
\begin{equation}
\cS (\Delta_{\Gamma}  + \cW)\cS [\sigma] = \cS [f] 
\end{equation}
%i.e. discretizing the integral operator as a composition of~$\cS$ followed by
%$\Delta_{\Gamma} + \cW$ followed by $\cS$,
results in a numerical null space proportional to the number of patches in the
triangulation and the order of discretization nodes on each triangle. The null
space is purely a numerical artifact and can be attributed to the particular
choice of discretization for representing the surface and the density. This
issue can be remedied by using a basis and discretization that enforces
smoothness of solutions across adjacent triangles, however, this introduces
additional complexity for discretizing the surface and representing the
solution. Alternatively, the issue can instead be remedied by solving the
mathematically equivalent integral equation
\begin{equation}
  \label{eq:ie-dis0}
\cS(\Delta_{\Gamma} + \cW)\cS [\sigma]  = \nabla \cdot \cS  \nabla_{\Gamma} \cS [\sigma] + \cS \cW\cS [\sigma] = \cS [f].
\end{equation}
In the above form, the operator on the left is discretized via the composition
of~$\cS$ followed by~$\nabla_{\Gamma}$, followed by an application of~$\cS$ on
the resulting vector density, and lastly followed by ($\nabla \cdot$). The
compact term~$\cS \cW\cS$ is the discretized analogously and added to the
previous computation. In the rest of the paper, the integral operators are
assumed to be discretized in the order they are written, from right-to-left.
We will denote the expanded operator above as
\begin{equation}
  \cA_1 = \nabla \cdot \cS  \nabla_{\Gamma} \cS  + \cS \cW\cS.
\end{equation}

For comparison, we also present results related to the 
numerical discretizations of the alternative integral representations for
solving this problem:
\begin{equation}
\label{eq:ie-dis1}
\cS (\Delta_{\Gamma} + \cW) \cS [\sigma]  = -\frac{\sigma}{4} + \mathcal{D}^2[\sigma] -
\mathcal{S}\left(\mathcal{S}'' + \mathcal{D}' +2H\mathcal{S}'
\right)[\sigma] = \cS [f] \, ,
\end{equation}
and
\begin{equation}
\label{eq:ie-dis2}
(\Delta_{\Gamma} + \cW) \cS^2 [\sigma] = -\frac{\sigma}{4} + (\cS')^2[\sigma] - (\cS''+\cD' + 2H\cS')\cS [\sigma] + \cW \cS^2 [\sigma]= f \, .
\end{equation}
The operators in these equations will be denoted as~$\cA_2$ and~$\cA_3$, where
\begin{equation}
  \begin{aligned}
    \cA_2 &= -\frac{1}{4}\cI + \mathcal{D}^2 -
    \mathcal{S}\left(\mathcal{S}'' + \mathcal{D}' +2H\mathcal{S}'
    \right), \\
    \cA_3 &= -\frac{1}{4}\cI + (\cS')^2 
    - (\cS''+\cD' + 2H\cS')\cS  + \cW \cS^2.
      \end{aligned}
\end{equation}

In the following, suppose that~$\tau_i$, for $i=1,2,\ldots N$, are estimates  of
the solution to the linear system on the current iteration of GMRES. We will
consider solving discretized versions of the integral equations
in~\eqref{eq:ie-dis0},~\eqref{eq:ie-dis1}, and~\eqref{eq:ie-dis2}.
Across all three integral equations, there are seven different layer
potential operators that appear: $\mathcal{S}$, $\mathcal{D}$,
$\mathcal{S}'$, $\partial_{1} \cS$, $\partial_{2} \cS$, $\partial_{3} \cS$,
and~$\mathcal{S}'' + \mathcal{D}'$, where $\partial_{j} \cS$ denotes a partial
derivative of
the single layer potential along the coordinate direction~$x_j$.

Referring to~\eqref{eq:nys-lcq2}, let $(\tilde{w}_{ij}^{\mathcal{K}})$ for
$\mathcal{K} \in \{ \mathcal{S},\mathcal{D}, \partial_{1} \cS, \partial_{2} \cS,
\partial_{3} \cS, \mathcal{S}',(\mathcal{S}''+\mathcal{D}') \}$ denote the
effective quadrature corrections for the corresponding discretized layer
potentials. In order to use a vector version of the FMM, it is essential that
the far part of the layer potentials be computed using the same set of source
locations. This is ensured by choosing oversampled quadrature nodes and
weights~$\{\bs_{i}, \wover_{i} \}_{i=1}^{\Nover}$ that are accurate for the far
field part of all seven kernels. 

For some of the kernels above such as $\cS', \cD, \cD'$, and $\cD''$, the normals need to be evaluated at the oversampled source locations. In order to evaluate $\bn(\bs_{j})$ for $\bs_{j} \in \Gamma_{\ell}$, this is done by using the interpolated values of  $\partial_{u} \bX^{\ell}$, and $\partial_{v} \bX^{\ell}$ as follows, 
\begin{equation}
\bn(\bs_{j}) = \frac{\partial_{u} \bX^{\ell}(u_{j},v_{j}) \times \partial_{v} \bX^{\ell}(u_{j},v_{j})}{|\partial_{u} \bX^{\ell}(u_{j},v_{j}) \times \partial_{v} \bX^{\ell}(u_{j},v_{j})|} \, ,
\end{equation}
where $(u_{j},v_{j}) \in T_{0}$ are the local coordinates of $\bs_{j}$ on $\Gamma_{\ell}$. The reason for using interpolated values of $\partial_{u} \bX^{\ell}$, and $\partial_{v} \bX^{\ell}$ instead of directly interpolating the normals is the smoother nature of the tangential derivatives as compared to the surface normal owing to the normalization factor $1/|\partial_{u} \bX \times \partial_{v} \bX|$ in the normal vector. In the \texttt{fmm3dbie} package, only first derivative information of the charts are stored for representing the surface. The mean curvature $H(\bx_{i})$ is obtained by computing second derivative of the charts  $\bX^{\ell}$ using spectral differentiation of $\partial_{u} \bX^{\ell}$, and $\partial_{v}  \bX^{\ell}$ for $\bx_{i} \in \Gamma_{\ell}$. 
Finally, for any function $f$ sampled on $\Gamma$, $\tilde{f}_{j}$, $j=1,2,\ldots \Nover$, denotes the samples of the function at the oversampled discretization nodes.

\subsection{FMM-accelerated application of~$\cA_1$}

In this section, we discuss the FMM accelerated evaluation of 
\begin{equation}
\cA_{1}[\tau] = \nabla \cdot \cS \nabla_{\Gamma} \cS [\tau] + \cS \cW\cS [\tau] \, .
\end{equation}
Let $\phi^{(1)}$ denote the sum
\begin{equation}
\phi^{(1)}(\bx) = \sum_{j=1}^{\Nover} \frac{\wover_{j} \tilde{\tau}_{j}}{|\bx-\bs_{j}|} \, .
\end{equation}
Using the FMM, we compute $\phi^{(1)}(\bx_{i})$ and $\nabla \phi^{(1)}(\bx_{i})$, $i=1,2\ldots N$, where $c_{1,j} = \wover_{j} \tilde{\tau}_{j}$, and $\bv_{1,j} = 0$, $j=1,2,\ldots \Nover$. 
Let $\nabla_{\ell} S$ and $S$ denote the discretized versions of $\partial_{\ell} \cS$ and $\cS$ respectively, where $\partial_{\ell} \cS$ as before is the gradient of $\cS$ along the coordinate direction $\ell$, with $\ell=1,2$, or $3$. Then
\begin{equation}
\begin{aligned}
S[\tau](\bx_{i}) &= \phi^{(1)}(\bx_{i}) \, \, + \,\, \sum_{\substack{j=1\\ \bx_{j} \in \textrm{Near}(\bx_{i})}}^{N} \tilde{w}_{ij}^{\mathcal{S}} \tau(\bx_{j}) \, , \\
\partial_{\ell} S[\tau](\bx_{i}) &= \partial_{\ell}\phi^{(1)}(\bx_{i}) \, \, + \,\, \sum_{\substack{j=1\\ \bx_{j} \in \textrm{Near}(\bx_{i})}}^{N} \tilde{w}_{ij}^{\partial_{\ell} \mathcal{S}} \tau(\bx_{j}) \, .
\end{aligned}
\end{equation}
Let $\nabla S = [\partial_{1} S; \partial_{2} S; \partial_{3} S]$ denote the discretized version of $\nabla \cS$, and
let $\bmu = [\mu^{(1)}; \mu^{(2)}; \mu^{(3)}]$ denote the discretized version of $\nabla_{\Gamma} \cS$, which can be expressed in terms of $\partial_{\ell} S$, and tangential derivatives $\partial_{u} \bX$, and $\partial_{v} \bX$ as follows
\begin{equation}
\begin{aligned}
\bmu(\bx_{i}) &= (\nabla S (\bx_{i}) \cdot  \partial_{u} \bX^{\ell}(u_{i},v_{i})) \left( g^{11} \partial_{u} \bX^{\ell}(u_{i},v_{i}) + g^{12} \partial_{v} \bX^{\ell}(u_{i},v_{i})  \right) + \\
& \quad (\nabla S (\bx_{i}) \cdot  \partial_{v} \bX^{\ell}(u_{i},v_{i})) \left( g^{21} \partial_{u} \bX^{\ell}(u_{i},v_{i}) + g^{22} \partial_{v} \bX^{\ell}(u_{i},v_{i})  \right) \, ,
\end{aligned}
\end{equation} 
where $\bx_{i} \in \Gamma_{\ell}$, $(u_{i},v_{i}) \in T_{0}$ are the local coordinates of $\bx_{i}$ on $\Gamma_{\ell}$, and $g^{ij}$, $i,j=1,2$, are the components of the inverse of the metric tensor which are computed using
$\partial_{u} \bX^{\ell}(u_{i},v_{i})$, and $\partial_{v} \bX^{\ell} (u_{i},v_{i})$. 

Given the samples of $\bmu$ at the oversampled nodes, let $\bzeta=[\zeta^{(1)}; \zeta^{(2)}; \zeta^{(3)}]$ denote the sum
 \begin{equation}
 \bzeta(\bx) = \sum_{j=1}^{\Nover} \frac{\wover_{j} \tilde{\bmu}_{j}}{|\bx_{i} - \bs_{j}|} \, .
 \end{equation}
 Using a second call to the vector version of the FMM, we compute $\bzeta(\bx_{i})$, and its gradient $\nabla \bzeta(\bx_{i})$, $i=1,2\ldots N$, where $c_{\ell,j} = \wover_{j} \zeta^{(\ell)}_{j}$, and $\bv_{\ell,j}=0$, $\ell=1,2,3$, and $j=1,2\ldots \Nover$. The discretized version of $\nabla \cdot \cS[\bmu]$ denoted by $\nabla \cdot S [\bmu]$ is given by
 \begin{equation}
 \nabla \cdot S [\bmu](\bx_{i}) = \sum_{\ell=1}^{3} \bigg(  \partial_{\ell} \zeta^{(\ell)}(\bx_{i}) + \sum_{\substack{j=1\\ \bx_{j} \in \textrm{Near}(\bx_{i})}}^{N} \tilde{w}_{ij}^{\partial_{\ell} \cS} \mu^{(\ell)}_{j} \bigg) \, .
 \end{equation}

In order to handle the inclusion of $\cS \cW \cS [\tau]$, we precompute and store the discretized version of $\cS[1]$ denoted by $\phi_{0} = S[1]$. The computation of $\phi_{0}$ is identical to the computation of $S[\tau]$. This cost of evaluating $\phi_{0}$ is then amortized across the number of GMRES iterations. Let $\eta$ be the constant given by
\begin{equation}
\eta = \frac{\sum_{j=1}^{N} S[\tau](\bx_{j}) w_{j}}{ \sum_{j=1}^{N} w_{j}} \, .
\end{equation}
Combining the computations above, the discretized version of the integral operator in~\cref{eq:ie-dis0} acting on the density $\tau$ is given by
\begin{equation}
\nabla \cdot S \nabla_{\Gamma} S[\tau] + SWS[\tau] = \nabla \cdot S[\bmu] + \eta \phi_{0} \, .
\end{equation}

\subsection{FMM-accelerated application of~$\cA_2$}

In this section, we discuss the FMM-accelerated evaluation of
\begin{equation}
\cA_{2} [\tau] = -\frac{\tau}{4} + \mathcal{D}^2[\tau] -
\mathcal{S}\left(\mathcal{S}'' + \mathcal{D}' +2H\mathcal{S}'
\right)[\tau] + \mathcal{SWS}[\tau] \, .
\end{equation}
Let $\phi^{(1)}, \phi^{(2)}$ denote the sums
\begin{equation}
\begin{aligned}
\phi^{(1)}(\bx) &= \sum_{j=1}^{\Nover} \frac{\wover_{j} \tilde{\tau}_{j}}{|\bx-\bs_{j}|} \, ,\\
\phi^{(2)}(\bx) &= -\sum_{j=1}^{\Nover} \wover_{j} \tilde{\tau}_{j} \bn(\bs_{j}) \cdot \nabla_{\bx} \frac{1}{|\bx-\bs_{j}|} \, ,
\end{aligned}
\end{equation}
where $\bn(\bs_{j})$, $j=1,2,\ldots \Nover$, are the interpolated values of the normal vector at the oversampled nodes.
Using the vector version of the FMM, we compute $\phi^{(\ell)}(\bx_{i})$, $\nabla \phi^{(\ell)}(\bx_{i})$, and $\nabla \nabla \phi^{(\ell)}(\bx_{i})$, $\ell=1,2$, $i=1,2\ldots N$, where $c_{1,j} = \wover_{j} \tilde{\tau}_{j}$, $c_{2,j} = 0$, $\bv_{1,j} = 0$, and $\bv_{2,j} = \wover_{j} \tilde{\tau}_{j} \bn(\bs_{j})$, $j=1,2,\ldots \Nover$.
Let $S,D,S', S'' + D'$ denote the discretized versions of $\mathcal{S}$, $\mathcal{D}$, $\mathcal{S}'$, and $\mathcal{S}'' + \mathcal{D}'$ respectively. Then,
\begin{equation}
  \begin{aligned}
S[\tau](\bx_{i}) &= \phi^{(1)}(\bx_{i}) \, \, + \,\, \sum_{\substack{j=1\\ \bx_{j} \in \textrm{Near}(\bx_{i})}}^{N} \tilde{w}_{ij}^{\mathcal{S}} \tau(\bx_{j}) \, , \\
D[\tau](\bx_{i}) &= \phi^{(2)}(\bx_{i})\,\, + \,\, \sum_{\substack{j=1\\ \bx_{j} \in \textrm{Near}(\bx_{i})}}^{N} \tilde{w}_{ij}^{\mathcal{D}} \tau(\bx_{j}) \, , \\
S'[\tau](\bx_{i}) &= \bn(\bx_{i}) \cdot \nabla \phi^{(1)}(\bx_{i}) \,\,  + \,\, \sum_{\substack{j=1\\ \bx_{j} \in \textrm{Near}(\bx_{i})}}^{N} \tilde{w}_{ij}^{\mathcal{S}'} \tau(\bx_{j}) \, ,
  \end{aligned}
\end{equation}
and
\begin{multline}    
(S'' + D')[\tau](\bx_{i}) =\\
 \bn(\bx_{i}) \cdot \nabla \nabla \phi^{(1)}(\bx_{i}) \cdot \bn(\bx_{i}) + \bn(\bx_{i}) \cdot \nabla \phi^{(2)}(\bx_{i}) \,\,+ \,\, \sum_{\substack{j=1\\ \bx_{j} \in \textrm{Near}(\bx_{i})}}^{N} \tilde{w}_{ij}^{(\cS'' + \cD')} \tau(\bx_{j}) \, .
\end{multline}

Given these four quantities, let $\mu^{(1)} = D[\tau]$, and let $\mu^{(2)}$ be given by
\begin{equation}
\mu^{(2)}(\bx_{i}) = -(S'' + D')[\tau](\bx_{i}) - 2H(\bx_{i}) S'[\tau](\bx_{i}) + \frac{\sum_{j=1}^{N} w_{j} S[\tau](\bx_{j})}{\sum_{j=1}^{N}  w_{j}} \, .
\end{equation}
The last term in the sum above is the discretized version of $\mathcal{WS}[\tau]$. 

Given the interpolated values of the densities $\mu^{(\ell)}$ at the oversampled nodes, and let $\zeta^{(\ell)}$, $\ell=1,2$, denote the sums
\begin{equation}
\begin{aligned}
\zeta^{(1)}(\bx) &= -\sum_{j=1}^{\Nover} \wover_{j} \tilde{\mu}^{(1)}_{j} \bn(\bs_{j}) \cdot \nabla_{\bx} \frac{1}{|\bx-\bs_{j}|} \, , \\
\zeta^{(2)}(\bx) &= \sum_{j=1}^{\Nover} \frac{\wover_{j} \tilde{\mu}^{(2)}_{j}}{|\bx-\bs_{j}|} \, .
\end{aligned}
\end{equation}
Using a second call to the vector version of the FMM, we compute $\zeta^{(\ell)}(\bx_{i})$, $\ell=1,2$, $i=1,2\ldots N$, where $c_{1,j} = 0$, $c_{2,j} = \wover_{j} \tilde{\mu}^{(1)}_{j}$, $\bv_{1,j} = \wover_{j} \tilde{\mu}^{(2)}_{j} \bn(\bs_{j})$ and $\bv_{2,j} = 0$, $j=1,2,\ldots \Nover$. As before, let $D[\mu^{(1)}]$, and $S[\mu^{(2)}]$ denote the discretized versions of $\cD[\mu^{(1)}]$ and $\cS[\mu^{(2)}]$, then 
\begin{equation}
\begin{aligned}
D[\mu^{(1)}](\bx_{i}) &= \zeta^{(1)}(\bx_{i}) \, \, + \,\, \sum_{\substack{j=1\\ \bx_{j} \in \textrm{Near}(\bx_{i})}}^{N} \tilde{w}_{ij}^{\mathcal{D}} \mu^{(1)}(\bx_{j}) \, , \\
S[\mu^{(2)}](\bx_{i}) &= \zeta^{(2)}(\bx_{i})\,\, + \,\, \sum_{\substack{j=1\\ \bx_{j} \in \textrm{Near}(\bx_{i})}}^{N} \tilde{w}_{ij}^{\mathcal{S}} \mu^{(2)}(\bx_{j}) \, .
\end{aligned}
\end{equation}
Finally, the discretized version of the integral operator in~\cref{eq:ie-dis1} acting on the density $\tau$, is given by
\begin{equation}
-\frac{\tau}{4} + D^2[\tau] -
S\left(S'' + D' +2HS'
\right)[\tau] +SWS[\tau] = -\frac{\tau}{4} + S[\mu^{(2)}] - D[\mu^{(1)}]  \, .
\end{equation}

\subsection{FMM-accelerated application of~$\cA_3$}

In this section, we discuss the FMM-accelerated evaluation of
\begin{equation}
\cA_{3} [\tau] = -\frac{\tau}{4} + (\cS')^2[\tau] - (\cS''+\cD' + 2H\cS')\cS [\tau] + \cW \cS^2 [\tau] \,.
\end{equation}
Let $\phi^{(1)}$, denote the sum
\begin{equation}
\begin{aligned}
\phi^{(1)}(\bx) &= \sum_{j=1}^{\Nover} \frac{\wover_{j} \tilde{\tau}_{j}}{|\bx-\bs_{j}|} \, .
\end{aligned}
\end{equation}
Using the FMM, we compute $\phi^{(1)}(x_{i})$, $\nabla \phi^{(1)}(x_{i})$ , $i=1,2\ldots N$, where $c_{1,j} = \wover_{j} \tilde{v}_{j}$, and $\bv_{1,j} = 0$, $j=1,2,\ldots \Nover$.
As before, let $S$, and $S'$ denote the discretized versions of $\mathcal{S}$, and $\mathcal{S}'$ respectively. Then
\begin{equation}
\begin{aligned}
S[\tau](\bx_{i}) &= \phi^{(1)}(\bx_{i}) \, \, + \,\, \sum_{\substack{j=1\\ \bx_{j} \in \textrm{Near}(\bx_{i})}}^{N} \tilde{w}_{ij}^{\mathcal{S}} \tau(\bx_{j}) \, , \\
S'[\tau](\bx_{i}) &= \bn(\bx_{i}) \cdot \nabla \phi^{(1)}(\bx_{i}) \,\,  + \,\, \sum_{\substack{j=1\\ \bx_{j} \in \textrm{Near}(\bx_{i})}}^{N} \tilde{w}_{ij}^{\mathcal{S}'} \tau(\bx_{j}) \, .
\end{aligned}
\end{equation}

Given these two quantities, let $\mu^{(1)} = S'[\tau]$, and $\mu^{(2)}=S[\tau]$. 
Let $\zeta^{(\ell)}$, $\ell=1,2,3$ denote the sums
\begin{equation}
\begin{aligned}
\zeta^{(1)}(\bx) &= \sum_{j=1}^{\Nover} \frac{\wover_{j} \tilde{\mu}^{(1)}_{j}}{|\bx-\bs_{j}|} \, ,\\
\zeta^{(2)}(\bx) &= \sum_{j=1}^{\Nover} \frac{\wover_{j} \tilde{\mu}^{(2)}_{j}}{|\bx-\bs_{j}|} \, , \\
\zeta^{(3)}(\bx) &= -\sum_{j=1}^{\Nover} \wover_{j} \tilde{\mu}^{(2)}_{j} \bn(\bs_{j}) \cdot \nabla_{\bx} \frac{1}{|\bx-\bs_{j}|} \, . 
\end{aligned}
\end{equation}
Using a second call to the vector version of the FMM, we compute $\zeta^{(\ell)}(x_{i})$, $\ell=1,2,3$, $i=1,2\ldots N$, where $c_{1,j} = \wover_{j} \tilde{\mu^{(1)}}_{j}$, $c_{2,j} = \wover_{j} \tilde{\mu}^{(2)}_{j}$, $c_{3,j} = 0$,  $\bv_{1,j} = \bv_{2,j} = 0$, and  $\bv_{3,j} = \wover_{j} \tilde{\mu}^{(2)}_{j} \bn(\bs_{j})$, $j=1,2,\ldots \Nover$. As before, let $S'[\mu^{(1)}]$, $S[\mu^{(2)}]$, $S'[\mu^{(2)}]$ and $(S'' + D')[\mu^{(2)}]$ denote the discretized versions of $\cS[\mu^{(1)}]$, $\cS[\mu^{(2)}]$, $\cS'[\mu^{(2)}]$ and $(\cS'' + \cD')[\mu^{(2)}]$ respectively. Then 
\begin{equation}
\begin{aligned}
S'[\mu^{(1)}](\bx_{i}) &= \bn(\bx_{i}) \cdot \nabla \zeta^{(1)}(\bx_{i}) \, \, + \,\, \sum_{\substack{j=1\\ \bx_{j} \in \textrm{Near}(\bx_{i})}}^{N} \tilde{w}_{ij}^{\mathcal{S'}} \mu^{(1)}(\bx_{j}) \, , \\
S[\mu^{(2)}](\bx_{i}) &= \zeta^{(2)}(\bx_{i})\,\, + \,\, \sum_{\substack{j=1\\ \bx_{j} \in \textrm{Near}(\bx_{i})}}^{N} \tilde{w}_{ij}^{\mathcal{S}} \mu^{(2)}(\bx_{j}) \, , \\
S'[\mu^{(2)}](\bx_{i}) &= \bn(\bx_{i}) \cdot \nabla \zeta^{(2)}(\bx_{i}) \, \, + \,\, \sum_{\substack{j=1\\ \bx_{j} \in \textrm{Near}(\bx_{i})}}^{N} \tilde{w}_{ij}^{\mathcal{S'}} \mu^{(2)}(\bx_{j}) \, ,
\end{aligned}
\end{equation}
and
\begin{multline}
(S'' + D')[\mu^{(2)}](\bx_{i}) =\\
 \bn(\bx_{i}) \cdot \nabla \nabla \zeta^{(2)}(\bx_{i}) \cdot \bn(\bx_{i}) + \bn(\bx_{i}) \cdot \nabla \zeta^{(3)}(\bx_{i}) \,\,+ \,\, \sum_{\substack{j=1\\ \bx_{j} \in \textrm{Near}(\bx_{i})}}^{N} \tilde{w}_{ij}^{(\cS'' + \cD')} \mu^{(2)}(\bx_{j}) \, .
\end{multline}

Having computed $S[\mu^{(2)}] = S^2[\tau]$, let $\eta$ be the discretized version of $\cW \cS^2 [\tau]$ given by
\begin{equation}
\eta = \frac{\sum_{j=1}^{N} S[\mu^{(2)}](\bx_{j}) w_{j} }{\sum_{j=1}^{N} w_{j}} \, .
\end{equation}
Finally, the discretized version of the integral operator in~\cref{eq:ie-dis2} acting on the density $\tau$, is given by
\begin{multline}
-\frac{\tau}{4} + (S')^2[\tau] -
\left(S'' + D' +2HS'
\right)S[\tau] +WS^2[\tau] \\
= -\frac{\tau}{4} + S'[\mu^{(1)}] - (S'' + D')[\mu^{(2)}] - 2HS'[\mu^{(2)}] + \eta  \, .
\end{multline}

\subsection{Summary of FMM-accelerations}

The rate limiting step in the solution of the Laplace Beltrami problem is the
evaluation of the quadrature corrections $\tilde{w}_{ij}^{\cK}$, and the FMM
sums for all three integral equations. The computational performance of the fast
multipole method depends on the distribution of sources and targets (which is
identical for all three integral equations), whether the FMM sums are evaluating
charge sums, or dipole sums, or charges and dipole sums, and whether just the
potential is requested, or potentials and gradients, or potentials, gradients
and Hessians. The numerical implementation of each of the three integral
representations requires two calls to the vector FMM, where the sum of the total
number of densities is four for all three representations. We summarize the
specifics of the FMM sums for all three integral equations in
the~\cref{tab:summary}.

\begin{table}[t!]
\centering
\caption{ Summary of FMMs used in discretizations
of~\cref{eq:ie-dis0,eq:ie-dis1,eq:ie-dis2}}
\begin{tabular}{ccccc} \toprule
     {Representation} &  & \# of densities & Interaction kernel  & Output   \\ \hline
%    \multirow{3}{*}{\cref{rep1}} & 48  & 480   & 3.1E-2 & 1.1E+12 \\    
%    & 192 & 1920 & 6.8E-3 & 4.6E+10 \\
%    & 768 & 7680 & 9.0E-4 & 3.0E+7 \\ \hline
        \multirow{2}{*}{$\cA_1$} & 1st FMM  & 1   & Charge & Potential and Gradient \\    
    & 2nd FMM & 3 & Charge & Potential and Gradient \\ \midrule
            \multirow{2}{*}{$\cA_2$} & 1st FMM  & 2   & Charge, and Dipole & Potential, Gradient, and Hessian \\    
    & 2nd FMM & 2 & Charge, and Dipole & Potential \\ \midrule
             \multirow{2}{*}{$\cA_3$} & 1st FMM  & 1   & Charge & Potential and Gradient \\    
    & 2nd FMM & 3 & Charge, and Dipole & Potential, Gradient, and Hessian \\ \bottomrule
\end{tabular}
 \label{tab:summary}
\end{table}

\section{Harmonic vector fields} 
\label{sec:harm-vec}

For a surface $\Gamma$ with unit normal $\bn$, a harmonic vector field $\bH$ is
a tangential vector field which satisfies 
\begin{equation}
\nabla_{\Gamma} \cdot \bH = 0 \, , \quad \nabla_{\Gamma} \cdot (\bn
\times \bH) = 0 \, .
\end{equation} 
On a genus $g$ surface, it is well known that the space of harmonic
vector fields is $2g$ dimensional~\cite{dai_2014, epstein-2013,
  epstein-2010}. Moreover, it follows from the definition that if
$\bH$ is a harmonic vector field, then $\bn \times \bH$ is also a
harmonic vector field which is linearly independent of $\bH$.  Any
tangential $C^{1}$ vector field $\bF$ along a smooth surface admits a
Hodge decomposition~\cite{epstein-2010, epstein-2013}:
\begin{equation}
\bF = \nabla_{\Gamma} \alpha + \bn \times \nabla_{\Gamma} \beta + \bH \, ,
\end{equation}
where $\alpha$ and $\beta$ are scalar functions defined on the surface
$\Gamma$ and $\bH$ is a harmonic vector field.  Given a smooth
tangential vector field $\bF$, one can use the existence of
this decomposition to
compute linearly independent harmonic vector fields $\bH$ by solving
the following Laplace-Beltrami equations and subsequently computing $\bH$:
\begin{equation}
\begin{aligned}
\Delta_{\Gamma} \alpha &= \nabla_{\Gamma} \cdot \bF \, ,\\
 \Delta_{\Gamma} \alpha &= -\nabla_{\Gamma} \cdot (\bn \times \bF),  \\
\bH &= \bF - \nabla_{\Gamma} \alpha - \bn\times \nabla_{\Gamma} \beta \, .
\end{aligned}
\end{equation}

A convenient way of choosing such smooth vector fields $\bF$ is to first define
a vector field~$\bV \in \bR^{3}$ where each component is a low-degree random
polynomial, and then set $\bF = \bn \times \bn \times \bV$. We use scaled
Legendre polynomials for this purpose. The maximum degree of Legendre
polynomials should be at least $\lceil (2g)^{1/3} \rceil$ so that dimension of
the space of polynomials is greater than the dimension of the space of the
vector fields. For the examples in this paper, we choose the degree to be
$\lceil (2g)^{1/3} \rceil + 3$. Moreover, suppose that $\Gamma$ is contained in
$(-L,L)^3$, then each component is defined to be a scaled tensor product
Legendre polynomial $P_{\ell}(x/2L) P_{m} (y/2L) P_{n} (z/2L)$, where $P_{n}(x)$
is the standard Legendre polynomials of degree $n$ on $[-1,1]$.  The extra
factor of $2$ in the choice of Legendre polynomials is to avoid requiring
additional degrees of freedom to represent $\bF$ on the surface due to large
gradients of the Legendre polynomials near the end points $-1$ and $1$.

With high probability, different random polynomial vector fields $\bF$ result in
different projections on the different elements of the $2g$-dimensional space of
vector fields. It is easy to show that if $\bH$ is a harmonic vector field then
$\bn \times \bH$ is also a harmonic vector field which is linearly independent
of $\bH$.  Thus, in practice one can construct a basis for the harmonic vector
field subspace by computing $\bH$ for $g$ different vector fields, and the
remaining $g$ vector fields of the basis are defined as $\bn \times \bH$.

In the next section, we present a battery of numerical results demonstrating the
accuracy and efficiency of our FMM-accelerated solver and its applications.

\section{Numerical examples}
\label{sec:num}

In this section, we provide several numerical examples demonstrating the
accuracy and computational efficiency of our solver for the Laplace-Beltrami
problem. The code was implemented in Fortran and compiled using the GNU Fortran
12.2 compiler. We use the point-based FMMs from the \texttt{FMM3D} package
(\url{https://github.com/flatironinstitute/FMM3D}), and the locally-corrected
quadratures from the \texttt{fmm3dbie} package
(\url{https://github.com/fastalgorithms/fmm3dbie}). All CPU timings in these
examples were obtained on a laptop with 8 cores of an Intel i9 2.4 GHz
processor. 

We compare the performance of the three different approaches to solve the Laplace-Beltrami problem. 
We summarize the different integral representations, and the corresponding integral equations being solved
in~\cref{tab:rep}.

\begin{table}[t!]
  \centering
  \caption{Summary of integral representations and the corresponding integral equations.}
  \begin{tabular}{ccc}
Operator & Representation & Integral equation \\ \toprule
$\mathcal{A}_{1}$ & $u = \cS[\sigma]$ & $\nabla \cdot \cS \nabla_{\Gamma}
\cS[\sigma] + \cS \cW \cS [\sigma]  = f$ \\ \midrule
$\cA_{2}$ & $u=\cS[\sigma ]$ & $-\frac{\sigma}{4} + \mathcal{D}^2[\sigma] -
\mathcal{S}\left(\mathcal{S}'' + \mathcal{D}' +2H\mathcal{S}'
\right)[\sigma] + \mathcal{SWS}[\sigma] = f$ \\ \midrule
$\cA_{3}$ & $u = \cS^2[\sigma]$ & $-\frac{\sigma}{4} + (\cS')^2[\sigma] - (\cS''+\cD' + 2H\cS')\cS [\sigma] + 
\cW \cS^2 [\sigma] = f$ \\
\bottomrule
\end{tabular}
\label{tab:rep}
\end{table}

In the following, let $p$ denote the discretization order, and $N$ denote the
total number of discretization points on the surface. Let $\varepsilon$ denote
the tolerance for evaluating the quadrature corrections, and the tolerance for
computing the FMM sums. Let $t_{q}$ denote the time taken for computing the
quadrature corrections, $t_{mv}$ denote the time taken for applying the
discretized integral operators $\cA_{j}$, and let $t_{s}$ denote the
time taken to solve the linear system. Finally, let $\niter$ denote the number
of iterations required for the relative residual in GMRES to drop below the
prescribed tolerance $\epsgmres$. We cap the maximum number of GMRES iterations
at~$100$.

\subsection{Speed and accuracy performance}

To illustrate the accuracy of our solvers, consider the Laplace-Beltrami problem
on a sphere with data given by a random $80$th order spherical harmonic
expansion,
\begin{equation}
f(\theta,\phi) = \text{Re} \left(\sum_{n=1}^{80} \sum_{m=0}^{n} c_{n,m} Y_{n,m} (\theta,\phi) \right) \, ,
\end{equation}
where $c_{n,m}$ are complex numbers with real and imaginary parts uniformly sampled in $(-0.5,0.5)$.
The analytic solution for this problem is given by
\begin{equation}
u(\theta,\phi) = \text{Re} \left( \sum_{n=1}^{80} \sum_{m=0}^{n} -\frac{c_{n,m} Y_{n,m} (\theta,\phi)}{n(n+1)} \right) \, .
\end{equation}
The sphere is discretized by a stereographic projection of a triangulated cube.
Let $\ucomp$ denote the solution computed using any of the three integral
equations. Let $\varepsilon_{s}$ denote the relative error in the solution $u$
scaled by the by the boundary data given by
\begin{equation}
\varepsilon_{s} = \frac{|\ucomp - u|_{{\mathbb L}^2(\Gamma)}}{|f|_{{\mathbb L}^{2}(\Gamma)}} \approx \sqrt{\frac{\sum_{j=1}^{N} |\ucomp(\theta_{j},\phi_{j}) - u(\theta_{j},\phi_{j})|^2 w_{j}}{ \sum_{j=1}^{N} |f(\theta_{j}, \phi_{j})|^2 w_{j}}} \, .
\end{equation}
Here $(\theta_{j},\phi_{j})$ are the discretization nodes, and $w_{j}$ are
the corresponding quadrature weights for integrating smooth functions.

We set $\varepsilon=\epsgmres$, but use different values for $\varepsilon$, depending on the value of $p$. In particular, we use the following combination of parameters: $(p, \varepsilon) =$ $(5, 5 \times 10^{-7})$, 
$(7, 5 \times 10^{-9})$, and $(9, 5\times 10^{-11})$. 
Let $\rho_{c}$ denote the computed order of convergence given by $\log_{2}(\varepsilon_{s}(N)/\varepsilon_{s}(4N))$.

Referring to~\cref{tab-comp}, we observe that the solution converges at
$O(h^{p+0.5})$. The fluctuating order of convergence for smaller values of $N$
can be attributed to the lack of resolution of the data $f$. Given the presence
of order $0$ kernels like $\cS'$ and $\cD$, one may have expected the
discretized layer potentials to converge at
$O(h^{p-1})$~\cite{greengard2020fmm,Atkinson95}. However, the added regularity
in the solution obtained by applying $\cS$ or $\cS^2$, depending on the integral
representation used, could account for the improved order of convergence. 

With regards to CPU-time performance, we observe that $t_{q}, t_{mv}$, and
$t_{s}$ scale linearly with $N$ (fluctuations in $t_{mv}$ for lower values of
$N$ can be attributed to the cache effects). Among the three representations,
the operator $\cA_{3}$ requires the smallest $t_{q}$ since it requires
quadrature corrections for three kernels, as opposed to four for the other two
representations. On the other hand, the $t_{mv}$ is smallest for $\cA_{1}$ since
it requires only charge interactions for both the FMM calls,
see~\cref{tab:summary}. 

Due to the well-conditioned nature of the problem, $\niter=O(1)$ for $\cA_{2}$
and $\cA_{3}$, and thus the solve time $t_{s}$ is dominated by $t_{mv}$.
Moreover, $\niter$ is independent of the refinement of the mesh for $\cA_{2}$
and $\cA_{3}$. The ratio of $t_{q}$ and $t_{mv}$ lies between 10 and 40 for most
configurations; this seems to indicate that when the iteration count exceeds
$40$, the quadrature time is no longer the dominant cost of solving the linear
system. This might not be the case in many practical problems due to the
well-conditioned nature of the integral equation. 

We now turn our attention to the results for $\cA_{1}$. It is mathematically
equivalent to $\cA_{2}$, but the GMRES iteration stalls and the iteration
terminates because of reaching the maximum iteration count. This results in the
largest solve time for $\cA_{1}$, despite having the smallest $t_{mv}$ per
iteration. This behavior can be attributed to the lack of \emph{numerical
second-kindness} discussed in greater detail in~\cref{sec:num-secondkindness}
below. 

To summarize, owing to its significantly smaller $t_{q}$, and small $\niter$, $\cA_{3}$ has the best numerical performance amongst $\cA_{1}, \cA_{2}$, and $\cA_{3}$ for solving the Laplace-Beltrami problem. 

\begin{remark}
The operator $\cA_{3}$ also has an added analytical advantage over the other two
representations: The map from $f \to u$ is a pseudo-differential operator of
order $-2$. By representing the solution as $u=\cS^2[\sigma]$, the map from
$\sigma \to u$ is also an pseudo-differential operator of order $-2$. Thus,
$\cA_{3}$ tends to be more analytically faithful to the spectral properties of
operators which contain compositions of solutions to the Laplace-Beltrami
problem as compared to $\cA_{1}$ or $\cA_{2}$. This situation naturally arises
when using the generalized-Debye formulation for solving Maxwell's equations
with perfect conductor or dielectric boundary
conditions~\cite{epstein-2010,gendeb2}, and also in modeling of type-I
superconductors~\cite{superconductor1}.
\end{remark}

\afterpage{
\begin{table}[t!]
\centering
\caption{ Timing and accuracy results on the unit sphere for integral equations corresponding to
$\cA_{1}, \cA_{2}$, and $\cA_{3}$, discretized with $p=5,7$, and $9$ order patches. $N$ is the total number of discretization points, $t_{q}$ is the time taken for computing quadrature corrections in seconds, $t_{mv}$ is the time taken to apply the FMM-accelerated discretized integral operators in seconds, $t_{s}$ is the time taken to solve the linear system using GMRES in seconds, $\niter$ is the number of GMRES iterations required for the solution to converge to the prescribed GMRES tolerance, $\varepsilon_{s}$ is the accuracy of the computed solution, relative to the norm of the data, and $\rho_{c}$ is the estimated order of convergence.}
\begin{tabular}{c|cccccccc} \midrule
     {Operator} & $p$ &  {$N$} & $t_{q}$ & $t_{mv}$ & $t_{s}$ & $\niter$ & {$\varepsilon_{s}$}  & $\rho_{c}$  \\ \midrule
     
     \multirow{9}{*}{$\cA_{1}$} & 5 & 11520   & 9.5 & 0.5 & 61.1 & 100  & $5.5 \times 10^{-3}$ &   \\    
    & 5 &  46080 & 37.6 & 1.8 & 248.2 & 100 & $3.3 \times 10^{-4}$ & 4.1 \\
    & 5 &  184320 & 171.8 & 7.8 & 1096.6 & 100 & $6.5 \times 10^{-6}$ & 5.7 \\ \cmidrule{2-9}
    
    & 7 &  21504 & 35.8 & 1.4 & 198.9 & 100 & $2.2 \times 10^{-3}$ & \\
        & 7 &  86016 & 157.2 & 7.0 & 947.6 & 100 & $6.2 \times 10^{-5}$ & 5.1 \\
            & 7 &  344064 & 630.2 & 31.1 & 3931.3 & 100 & $2.2 \times 10^{-7}$ & 8.1 \\ \cmidrule{2-9}
    
    & 9 &  34560 & 206.0 & 6.0 & 976.9 & 100 & $1.1 \times 10^{-3}$ & \\
        & 9 &  138240 & 837.1 & 33.0 & 4001.8 & 100 & $7.4 \times 10^{-6}$ & 7.21 \\
            & 9 &  552960 & 2988.2 & 75.5 & 12517.0 & 100 & $6.8 \times 10^{-9}$ & 10.0 \\ \midrule

        \multirow{9}{*}{$\cA_{2}$} & 5 & 11520   & 6.8 & 1.4 & 16.8 & 5  & $1.4 \times 10^{-3}$ &   \\    
    & 5 &  46080 & 31.2 & 4.5 & 62.2 & 5 & $1.4 \times 10^{-4}$ & 3.3 \\
    & 5 &  184320 & 127.4 & 14.0 & 254.3 & 5 & $1.3 \times 10^{-7}$ & 10.0 \\ \cmidrule{2-9}
    
    & 7 &  21504 & 28.2 & 3.1 & 54.6 & 6 & $3.7 \times 10^{-4}$ & \\
        & 7 &  86016 & 116.6 & 15.7 & 228.5 & 6 & $5.9 \times 10^{-7}$ & 9.29 \\
            & 7 &  344064 & 499.1 & 47.8 & 850.0 & 6 & $3.0 \times 10^{-9}$ & 7.61 \\ \cmidrule{2-9}
    
    & 9 &  34560 & 156.5 & 8.9 & 245.8 & 7 & $3.9 \times 10^{-5}$ & \\
        & 9 &  138240 & 680.2 & 39.6 & 1046.7 & 7 & $5.4 \times 10^{-8}$ & 9.49 \\
            & 9 &  552960 & 2404.4 & 109.0 & 3429.2 & 7 & $7.4 \times 10^{-11}$ & 9.51 \\ \midrule

                \multirow{9}{*}{$\cA_{3}$} & 5 & 11520   & 6.0 & 0.9 & 11.8 & 5  & $1.4 \times 10^{-3}$ &   \\    
    & 5 &  46080 & 24.4 & 3.4 & 50.4 & 5 & $1.4 \times 10^{-4}$ & 3.3 \\
    & 5 &  184320 & 114.6 & 12.4 & 187.3 & 4 & $4.6 \times 10^{-7}$ & 8.2 \\ \cmidrule{2-9}
    
    & 7 &  21504 & 23.4 & 2.8 & 41.8 & 6 & $3.7 \times 10^{-4}$ & \\
        & 7 &  86016 & 96.0 & 13.6 & 199.4 & 6  & $5.9 \times 10^{-7}$ & 9.29 \\
            & 7 &  344064 & 377.4 & 57.1 & 722.3 & 5 & $4.0 \times 10^{-9}$ & 7.2 \\ \cmidrule{2-9}
    
    & 9 &  34560 & 123.8 & 8.5 & 196.1 & 7 & $3.9 \times 10^{-5}$ & \\
        & 9 &  138240 & 501.3 & 30.7 & 822.8 & 7 & $5.4 \times 10^{-8}$ & 9.49 \\
            & 9 &  552960 & 1755.4 & 104.3 & 2749.3 & 7 & $7.4 \times 10^{-11}$ & 9.51 \\ \midrule

\end{tabular}
  \label{tab-comp}
\end{table}

\clearpage
}

\subsection{Numerical second-kindness}
\label{sec:num-secondkindness}

\begin{figure}[t!]
  \centering
  \includegraphics[width=0.9\linewidth]{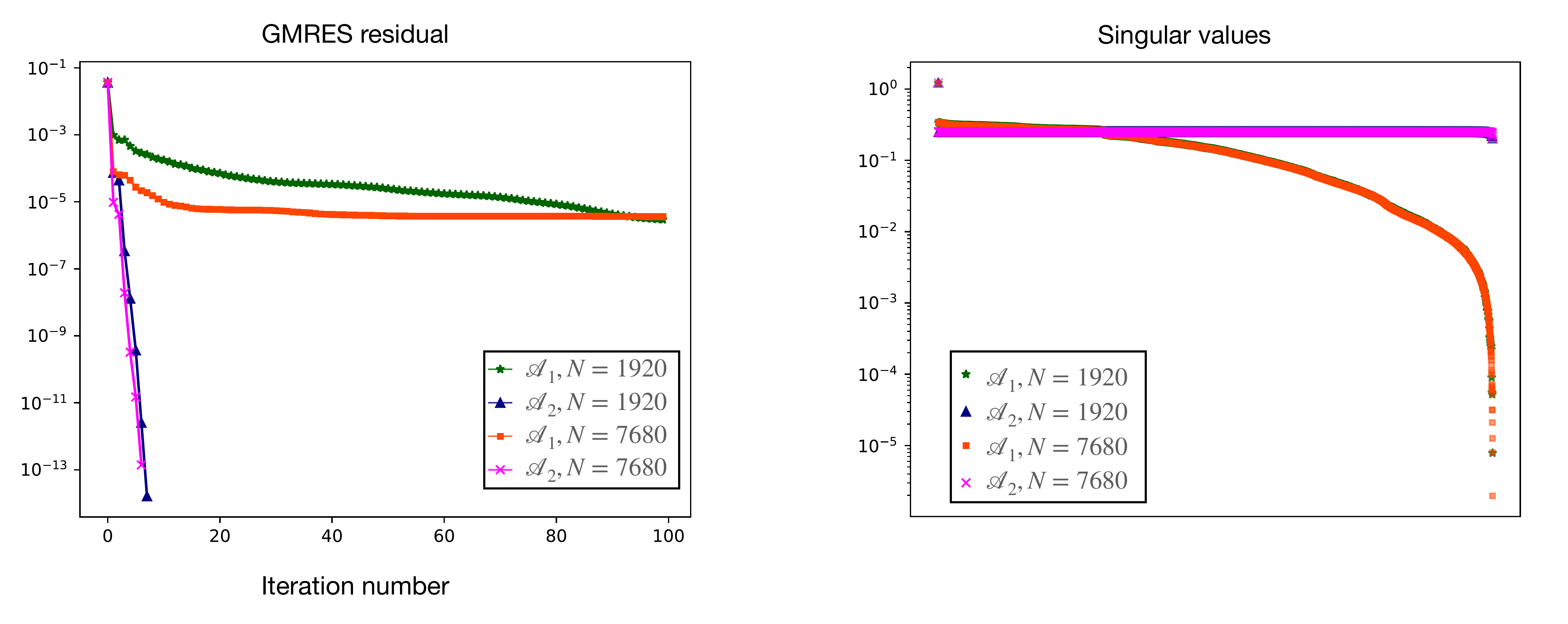}
  \caption{The singular values of the matrices for $\cA_{1}$, and $\cA_{2}$, for a sphere discretized with fourth order patches, and 
    $N=1920$, and $7680$ (left) and the GMRES residual as a
    function of iteration number (right). }
  \label{fig:num-secondkindness}
  \end{figure}
  
The poorer performance of $\cA_{1}$ in an iterative solver as compared to
$\cA_{2}$ can be attributed to the lack of \emph{numerical second-kindness} in
its implementation. A discretized numerical operator is defined to be
numerically second-kind if the identity term of the discretized integral
operator is handled explicitly as is the case for $\cA_{2}$ and $\cA_{3}$.
\emph{Any} numerical discretization of integral operators is inherently compact
since the operator is approximated via a finite-rank operator. This results in
the spectrum of the discretized integral operator clustering closely to the
origin which can  result in numerical conditioning issues. Explicitly adding the
identity term, as is the case in $\cA_{2}$, results in the spectrum of the
discretized integral operator clustering at $-1/4$. This \emph{greatly} enhances
the convergence of iterative methods such as GMRES~\cite{saad-1986}.

To avoid any additional errors due to using the FMM, we illustrate this behavior
on a smaller problem which would allow for using dense linear algebra routines.
Consider the Laplace-Beltrami operator on the sphere with data given by a random
fourth-order spherical harmonic expansion,
\begin{equation}
f(\theta,\phi) = \text{Re} \left(\sum_{n=1}^{4} \sum_{m=0}^{n} c_{n,m} Y_{n,m} (\theta,\phi) \right) \, ,
\end{equation}
where $c_{n,m}$ as before are complex numbers with real and imaginary parts
uniformly sampled in $(-0.5,0.5)$. We restrict our attention to $p=4$ and
$N=1920$ or $7680$. The linear system is then solved using GMRES with $\epsgmres
= 10^{-14}$, and $\varepsilon = 5\times 10^{-8}$. The matrices corresponding to
the integral equations are precomputed and at every iteration the matrix vector
product is computed using BLAS routines. Referring
to~\cref{fig:num-secondkindness}, we see that the GMRES residual for $\cA_{2}$
drops exponentially with iteration number, while that for $\cA_{1}$ converges at
a slower rate and eventually stalls. The plot of the singular values for the two
operators validates the lack of numerical second-kindness of $\cA_{1}$ as well
as the 
stability of condition number under refinement for the numerically second-kind
$\cA_{2}$. Finally, note that the smallest singular value for $\cA_{1}$
decreases approximately by a factor of $4$ upon refinement, indicating that it
might be proportional to area of the triangulated patches.

\subsection{Harmonic fields}

In this section, we present the numerical results for the computation of the
harmonic vector fields on two geometries: a twisted torus geometry whose boundary
$\Gamma$ is parametrized by $\bX: [0,2\pi]^2 \to \Gamma$ with
\begin{equation}
\bX(u,v) = \sum_{i=-1}^{2} \sum_{j=-1}^{2} \delta_{i,j} 
\begin{pmatrix}
\cos{v} \cos{((1-i)u + jv)} \\
\sin{v} \cos{((1-i)u + jv)} \\
\sin{((1-i)u + jv)} \\
\end{pmatrix} \, ,
\end{equation}
where the non-zero coefficients are
\begin{equation}
  \begin{aligned}
    \delta_{-1,-1} &=0.17,  &\qquad \delta_{-1,0} &= 0.11,  &\qquad
    \delta_{0,0} &= 1,\\
   \delta_{1,0} &= 4.5, & \delta_{2,0} &= -0.25, & \delta_{0,1} &=
   0.01,\\
   & & \delta_{2,1} &= -0.45 , & & 
  \end{aligned}
\end{equation}
and a genus-10 surface  obtained by first constructing the surface as a union of
rectangular faces.  A high-order mesh for a smooth version of the geometry is
obtained using the surface smoother algorithm of~\cite{vico2020surface},
see~\cref{fig:g10} for a graphical depiction. For both examples, the surface
is discretized with order $p=9$ patches. The quadrature and GMRES tolerances
were set to $\varepsilon = 5 \times 10^{-8}$, and $\epsgmres = 5 \times
10^{-10}$, respectively. In a slight abuse of notation, let $t_{s}$ now denote
the time taken for computing the Hodge decomposition of a given vector field,
which corresponds to solving two Laplace-Beltrami problems. Since the
computation of the quadrature corrections tend to be the dominant cost, we reuse
the computed quadrature corrections across the two solves. The solution to the
Laplace-Beltrami problems were computed using the representation $\cA_{3}$.

\begin{figure}[t!]
  \centering
\begin{subfigure}{.45\textwidth}
  \centering
  \includegraphics[width=.8\linewidth]{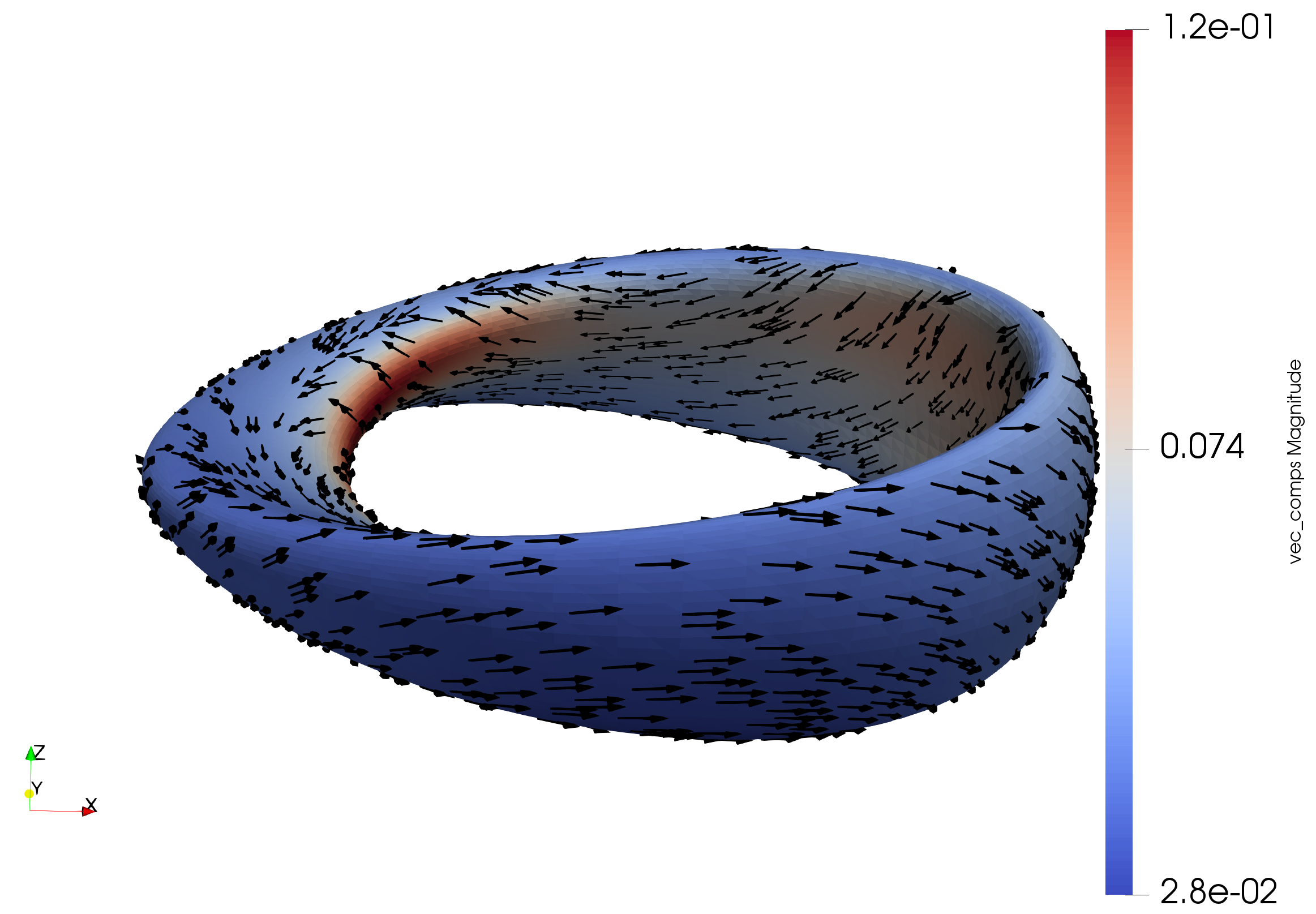}
  \caption{Tangential vector field, $\mathbf{F}$}
  \label{fig:htor1}
\end{subfigure}%
\quad 
\begin{subfigure}{0.45\textwidth}
  \centering
  \includegraphics[width=.8\linewidth]{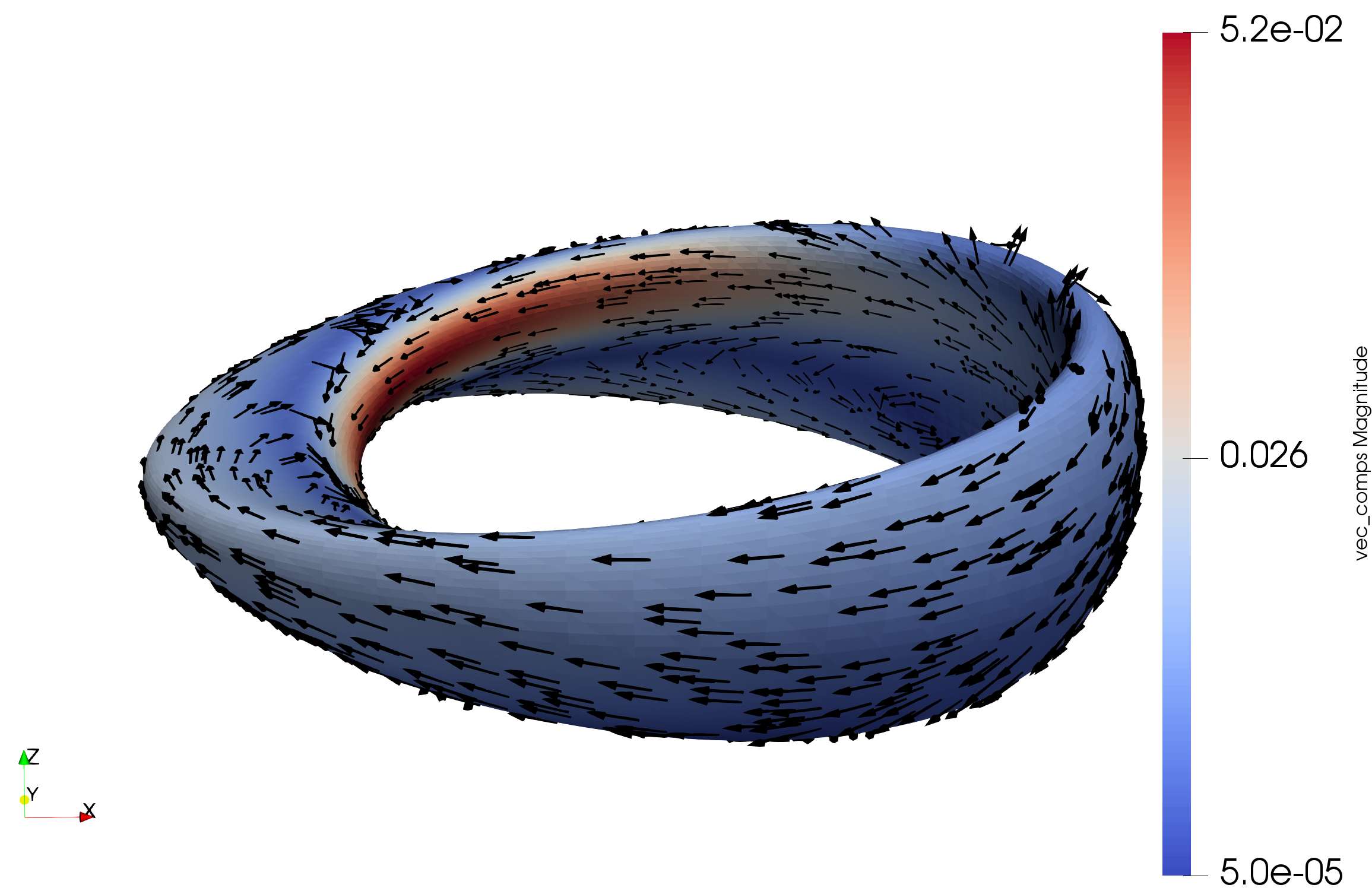}
  \caption{Curl-free component, $ \nabla_{\Gamma}\alpha$}
  \label{fig:htor2}
\end{subfigure}
\begin{subfigure}{.45\textwidth}
  \centering
  \includegraphics[width=.8\linewidth]{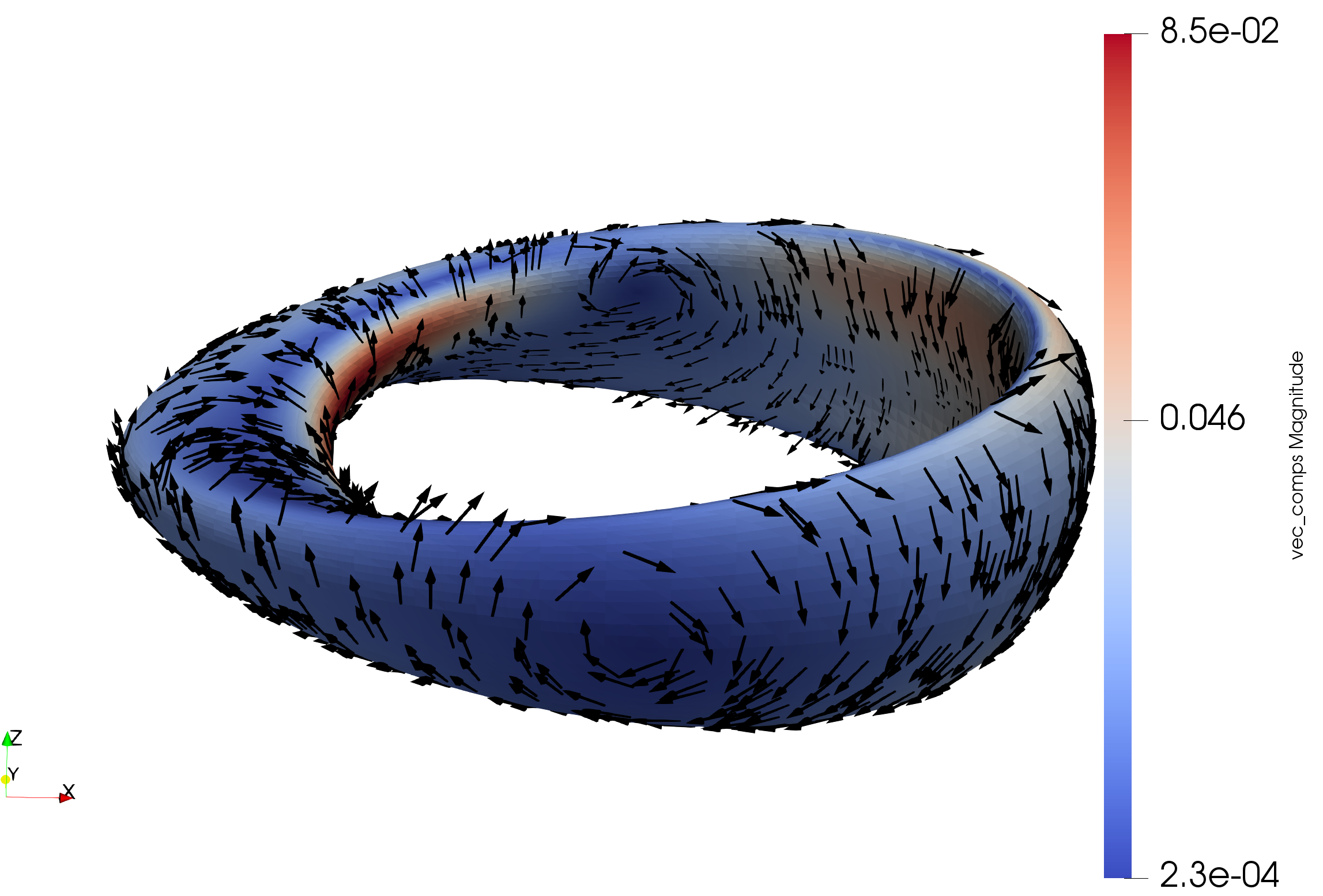}
  \caption{Div-free component, $\mathbf{n} \times \nabla_{\Gamma}\beta$  }
  \label{fig:htor3}
\end{subfigure}
\quad
\begin{subfigure}{.45\textwidth}
  \centering
  \includegraphics[width=.8\linewidth]{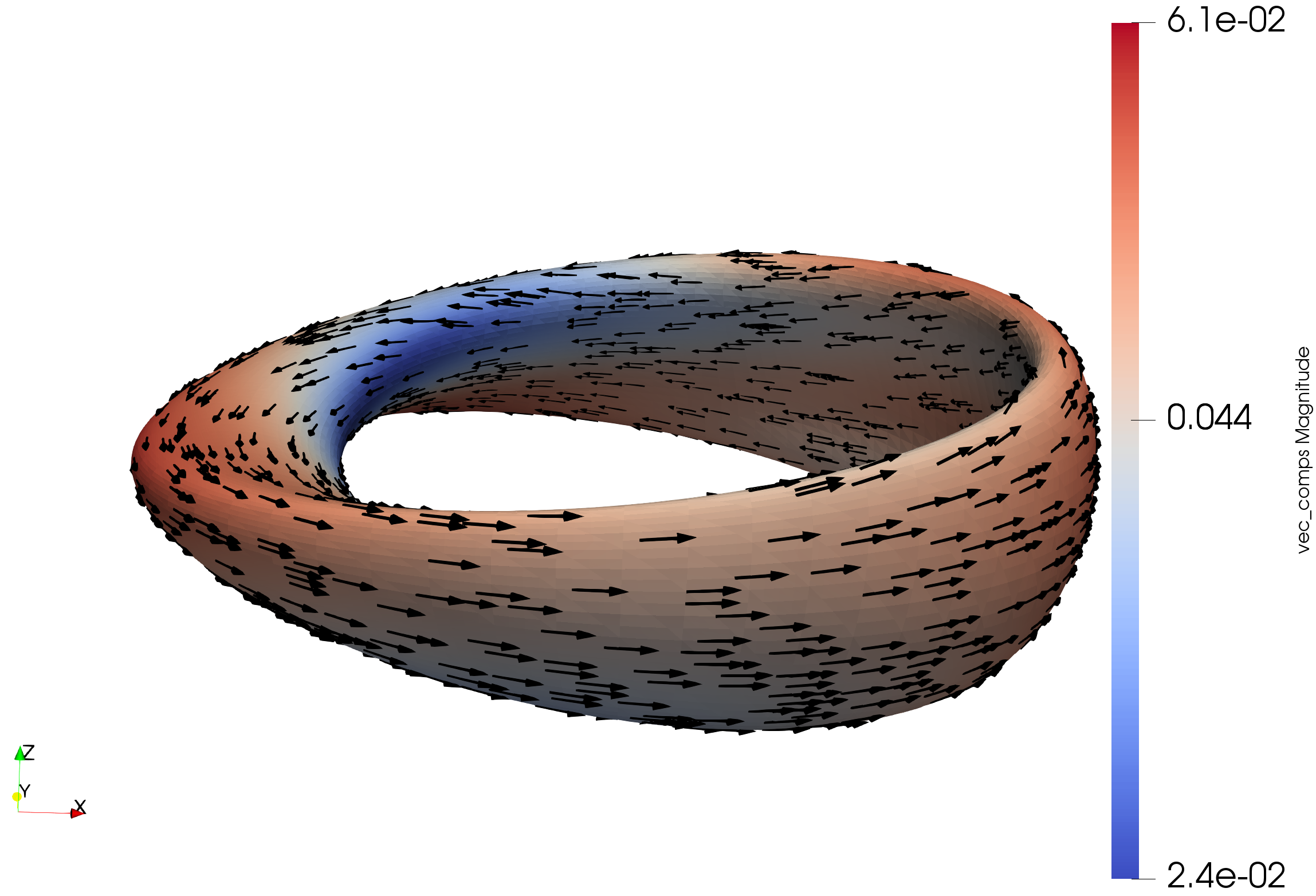}
  \caption{Harmonic component, $\mathbf{H}$}
  \label{fig:htor4}
\end{subfigure}
\caption{Hodge decomposition on a twisted torus with initial field
  $\bV$ given by a point current source at $\bx_{0}$.
  %: (a)
  %$\mathbf{F}$, the tangential projection of $\mathbf{V}$ (b)
  %curl-free component, $ \nabla_{\Gamma}\alpha$, (c) divergence free
  %component, $\mathbf{n} \times \nabla_{\Gamma}\beta$, (d) harmonic
  %field $\mathbf{H}$.
  The color bar represents the magnitude of the vector field and the
  arrows represent the direction of the field.}
\label{fig:htor}
\end{figure}

\begin{table}[b!]
  \centering
\begin{tabular}{cccccc} \toprule
    {$p$} & {$N$} & $\niter$ & $t_{s}$ & {$||\nabla_{\Gamma}\cdot \mathbf{\bH}||_{\mathbb{L}^{2}(\Gamma)}$  } &  {$||\nabla_{\Gamma}\cdot (\mathbf{n} \times \mathbf{\bH})||_{\mathbb{L}^{2}(\Gamma)}$}  \\ \midrule
    9  & 27,000 & (27,28) & 317.1 &$1.4 \times  10^{-2}$  & $1.0 \times 10^{-2}$     \\   
    9  & 108,000 & (27,28) & 1182.0 & $2.8 \times 10^{-4}$ & $1.5 \times 10^{-4}$   \\ \bottomrule
\end{tabular}
\caption{Results summary for computing harmonic vector fields on the twisted torus geometry. } \label{tab:hodgestell}
\end{table}

For the twisted torus example, let $\bV = \boldsymbol{\ell} \times (\bx -
\bx_{0})/| \bx -\bx_{0}|^{3}$, with $\bx_{0} = (0.2,0.2,0.2)$ and
$\boldsymbol{\ell} = (0,1,1)$, and then define $\bF = -\bn \times \bn \times
\bV$.  It is straightforward to see that a current source defined in the
exterior projects onto the harmonic vector fields for the torus.
In~\cref{fig:htor}, we plot the tangential vector field $\bF$, its curl-free
component, its divergence-free component and the harmonic vector field $\bH$.
The second linearly independent harmonic vector field can be obtained by
computing $\bn \times \bH$. In~\cref{tab:hodgestell}, we tabulate the number of
iterations for computing the curl-free and divergence-free parts of
decomposition denoted by $\niter$, and the solve time $t_{s}$. The error in the
computation of the harmonic vector fields is estimated by computing the
$\mathbb{L}^{2}$ norm of $\nabla_{\Gamma} \cdot \bH$ and $\nabla_{\Gamma} \cdot
(\bn \times \bH)$ using spectral differentiation. For this example $\| \bH \|
\approx 0.62498$, so the above error estimate serves as a proxy for the relative
error as well.

We also compute the 20 linearly independent harmonic vector fields on a genus 10
surface. For this geometry, we choose a random surface vector field $\bV$ whose
components are random tensor product Legendre polynomials as discussed
in~\cref{sec:harm-vec}. In~\cref{fig:g10}, we plot the surface along with these
20 linearly independent harmonic vector fields computed using our solver.  As
before, we verify the computations by calculating the $\mathbb{L}^{2}$ norms of
the surface divergence of $\bH$ and $\bn \times \bH$, where $\bH$ is the
computed harmonic field. In~\cref{tab:hodgeg10}, we tabulate the number of
iterations for computing the curl-free and divergence-free parts of
decomposition for one of the solves on the genus 10 surface, and also the
computation time for obtaining the hodge decomposition.

\begin{table}[b!]
  \centering
  \caption{Results summary for computing harmonic vector fields on the genus 10 geometry. }
    \begin{tabular}{cccccc} \toprule
    {$p$} & {$N$} & $\niter$ & $t_{s}$ & {$||\nabla_{\Gamma}\cdot \mathbf{H}||_{\mathbb{L}^{2}(\Gamma)}$  } &  {$||\nabla_{\Gamma}\cdot \mathbf{n} \times \mathbf{H}||_{\mathbb{L}^{2}(\Gamma)}$}  \\ \midrule
    9  & 388,800 & (19,19) & 2910.5 & $3.2 \times 10^{-2}$ & $5.6 \times 10^{-3}$      \\
    9  & 1,555,200 & (19,18) & 10298.0 & $5.8 \times 10^{-4}$ & $6.5 \times 10^{-5}$      \\\bottomrule
\end{tabular}
 \label{tab:hodgeg10}
\end{table}

The summarized results for both examples illustrate the linear scaling of the algorithm, and also the constancy of number of GMRES iterations under mesh refinement.

\afterpage{
\begin{figure}
\centering
\includegraphics[width=0.95\linewidth]{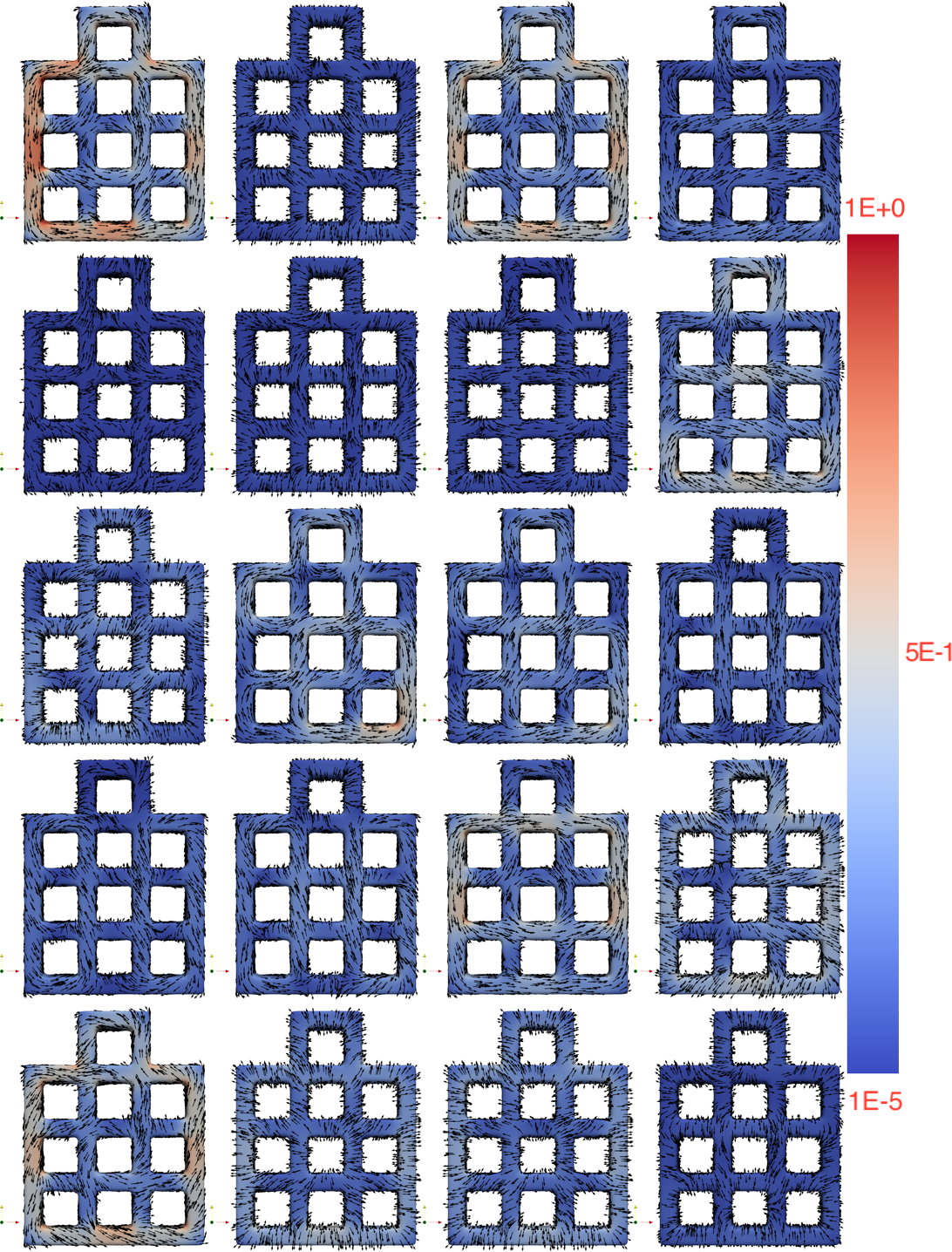}
    \caption{The 20 harmonic fields computed on a genus 10 geometry
    discretized with $N = 388800$ points. The
      color bar represents the magnitude of the vector field and the
      arrows represent the direction of the field.}
    \label{fig:g10}
\end{figure}

\clearpage
}

\section{Conclusions and future work}
\label{sec:conclusions}

In this paper, we have presented a high-order FMM-accelerated
iterative solver for the numerical solution of the Laplace-Beltrami
equation on complex smooth surfaces in three dimensions. The
Laplace-Beltrami problem~$\Delta_\Gamma \psi = f$ is converted to a
second-kind integral equation by representing the solution as
$\psi=\cS[\sigma]$ or $\psi=\cS^2[\sigma]$, and using appropriate Cald\'eron
identities. The resulting integral equation, which requires the
application of various Laplace layer potentials, is then discretized
using a high-order method with locally corrected quadratures for
computing the weakly-singular layer potentials; their application is
accelerated using fast multipole methods.

While the integral equation can be written in several equivalent
analytic forms, we demonstrate the necessity of maintaining numerical
second-kindness by explicitly isolating the identity term of the
Fredholm equation in order to avoid stagnation of iterative solvers
such as GMRES. We also illustrate that the using the representation
$\psi = \cS^2[\sigma]$ has the better numerical performance over
using $\psi = \cS[\sigma]$ and preconditioning with $\cS$.

Finally, we also presented numerical examples demonstrating the
computation of harmonic vector fields on surfaces using the
Laplace-Beltrami solver. These vector fields are both of pure
mathematical interest and are also required for solving problems in
type I superconductivity and electromagnetism on topologically
non-trivial geometries.

There are still several open questions that remain to be
addressed. These include coupling high-order discretizations of the
Laplace-Beltrami integral equations to fast direct solvers;
subsequently using these solvers to develop iterative solvers for the
solution of Maxwell's equations using the generalized Debye
formulation; and extending these ideas for the solution of surface
diffusion problems that arise in pattern formation and cell
biology. These are all ongoing areas of research.

\section{Data availability}
The code used in the numerical examples in this manuscript is available in two
publicly accessible git repositories:
\begin{itemize}
    \item \url{https://github.com/flatironinstitute/FMM3D}
    \item \url{https://github.com/fastalgorithms/fmm3dbie}
\end{itemize}

\section{Funding}
Dhwanit Agarwal's research was supported in part by UT-Austin CoLab and UTEN
Partnership, and the Portuguese Science and Technology Foundation under funding
~\#201801976/UTA18-001217. Michael O'Neil's research was  supported in part by
the Office of Naval Research under award numbers~\#N00014-17-1-2451
and~\#N00014-18-1-2307, and the Simons Foundation/SFARI (560651, AB).

\section{Conflict of interest}
The authors do not have any conflicts of interest.

\bibliographystyle{abbrv}
\bibliography{master}

\end{document}